\documentclass[12pt,a4paper,reqno]{amsart}
\usepackage{graphics,epic}
\usepackage{amsmath,amssymb, amsthm}
\usepackage[all,2cell]{xy}

\newtheorem{theorem}{Theorem}[section]
\newtheorem*{theorem*}{Theorem}

\newtheorem{lemma}[theorem]{Lemma}
\newtheorem{proposition}[theorem]{Proposition}
\newtheorem{corollary}[theorem]{Corollary}

\newtheorem*{conjecture*}{Conjecture}

\newtheorem{remark}[theorem]{Remark}

\newcommand{\n}{\mathfrak{N}}

\renewcommand{\hat}[1]{\widehat{#1}}

\newcommand{\id}{{\rm id}}

\newcommand{\Hom}{{\rm Hom}\,}

\newcommand{\End}{{\rm End}\,}

\newcommand{\Z}{\mathbb{Z}}

\newcommand{\C}{\mathbb{C}}

\def\C{{\mathbb C}}

\def\Z{{\mathbb Z}}

\def\1{{\bf 1}}

\def \End{{\rm End}}
\def \Hom{{\rm Hom}}

\def \Ind{{\rm Ind}}

\def \pf{\noindent {\bf Proof: \,}}
\def\theequation{5.\arabic{equation}}
\def \h{\mathfrak{h}}
\def \w{\omega}
\def \g{\mathfrak{g}}

\begin{document}

\title[Rationality and $C_2$-cofiniteness of certain diagonal coset VOAs]{Rationality and $C_2$-cofiniteness of certain diagonal coset vertex operator algebras}

%\vspace{0.5cm}
%\sffamily
\author{Xingjun Lin}
\address{Xingjun Lin,  School of Mathematics and Statistics, Wuhan University, Wuhan 430072, China.}
\thanks{X. Lin was supported by China NSF grant
11801419 and the starting research fund from Wuhan University (No. 413000076)}
\email{linxingjun88@126.com}
\begin{abstract}
In this paper, it is shown that the diagonal coset vertex operator algebra $C(L_{\mathfrak{g}}(k+2,0),L_{\mathfrak{g}}(k,0)\otimes L_{\mathfrak{g}}(2,0))$ is rational and $C_2$-cofinite in case  $\mathfrak{g}=so(2n), n\geq 3$ and $k$ is an admissible number for $\hat{\mathfrak{g}}$.  It is also shown that the diagonal coset vertex operator algebra $C(L_{sl_2}(k+4,0),L_{sl_2}(k,0)\otimes L_{sl_2}(4,0))$ is rational and $C_2$-cofinite in case $k$ is an admissible number for $\hat{sl_2}$. Furthermore, irreducible modules of $C(L_{sl_2}(k+4,0),L_{sl_2}(k,0)\otimes L_{sl_2}(4,0))$ are classified  in case  $k$ is a positive odd integer.
\end{abstract}
\maketitle
\section{Introduction}\label{intro}
\def\theequation{1.\arabic{equation}}
\setcounter{equation}{0}
Let $\g$ be a finite dimensional simple Lie algebra and $k$ be an admissible number for the affine Lie algebra $\hat{\g}$ \cite{KW1}. Then the irreducible quotient $L_{\g}(k, 0)$ of the vacuum module $V_{\g}(k, 0)$ of $\hat \g$ is a vertex operator algebra \cite{FZ}. For a positive integer $l$, the tensor product vertex operator algebra $L_{\g}(k,0)\otimes L_{\g}(l,0)$ is a $\hat \g$-module  of level $k+l$ with the diagonal action of $\hat \g$. We use
$C(L_{\g}(k+l,0),L_{\g}(k,0)\otimes L_{\g}(l,0))$ to denote the multiplicity space of $L_{\g}(k+l,0)$ in $L_{\g}(k,0)\otimes L_{\g}(l,0)$.
Then $C(L_{\g}(k+l,0),L_{\g}(k,0)\otimes L_{\g}(l,0))$ is a vertex operator algebra, which is called the diagonal coset vertex operator algebra.

 It is a longstanding conjecture that $C(L_{\g}(k+l,0),L_{\g}(k,0)\otimes L_{\g}(l,0))$ is rational and $C_2$-cofinite if $k$ is an admissible number for $\hat{\g}$ and $l$ is a positive integer (see Conjecture 7.2 of \cite{CL}). When $k=l=1$, rationality and $C_2$-cofiniteness of $C(L_{\g}(k+l,0),L_{\g}(k,0)\otimes L_{\g}(l,0))$ have been established for the case $\g=sl_{n}$ \cite{ACL}, \cite{JL1}, \cite{JL2} and the case $\g=so(n)$ \cite{AP}. In general case that $k$ is an admissible number for $\hat{\g}$ and $l=1$, it has been proved in \cite{ACL} that $C(L_{\g}(k+l,0),L_{\g}(k,0)\otimes L_{\g}(l,0))$ is rational and $C_2$-cofinite if $\g$ is simply-laced. It has also been proved in \cite{CL} that $C(L_{\g}(k-\frac{1}{2},0),L_{\g}(k,0)\otimes L_{\g}(-\frac{1}{2},0))$ is rational and $C_2$-cofinite  when $\g=sp_{2n}$ and $k$ is a positive integer. However, for $l\geq 2$, the conjecture is only confirmed in the special case $\g=sl_2$, $l=2$ \cite{A}, \cite{CFL} and  the special case $\g=E_8$, $l=2$ \cite{CL}, \cite{Lin}.

 In this paper, we study rationality and $C_2$-cofiniteness of  vertex operator algebras $C(L_{so(n)}(k+2,0),L_{so(n)}(k,0)\otimes L_{so(n)}(2,0))$ and $C(L_{sl_2}(k+4,0),L_{sl_2}(k,0)\otimes L_{sl_2}(4,0))$, which are closely related to Fateev-Zamolodchikov  $D_N$-parafermion algebras \cite{FaZa}, \cite{GS}. Explicitly, we show that $$ C(L_{so(2n)}(k+2,0),L_{so(2n)}(k,0)\otimes L_{so(2n)}(2,0))$$ is rational and $C_2$-cofinite if $n\geq 3$ and $k$ is an admissible number for $\hat{so(2n)}$ (see Theorem \ref{rationalD}). In the special case $n=3$, this implies that $$C(L_{sl_4}(k+2,0),L_{sl_4}(k,0)\otimes L_{sl_4}(2,0))$$ is rational and $C_2$-cofinite if  $k$ is an admissible number for $\hat{sl_4}$ (see Corollary \ref{rationala3}). We also show that $$C(L_{sl_2}(k+4,0),L_{sl_2}(k,0)\otimes L_{sl_2}(4,0))$$ is rational and $C_2$-cofinite if $ k$ is an admissible number for $\hat{sl_2}$ (see Theorem \ref{rationalmain}). Rationality and $C_2$-cofiniteness of $C(L_{so(2n+1)}(k+2,0),L_{so(2n+1)}(k,0)\otimes L_{so(2n+1)}(2,0))$ are also studied under certain assumption. More precisely, for a positive integer  $k$,  it is proved that $$C(L_{so(2n+1)}(k+2,0),L_{so(2n+1)}(k,0)\otimes L_{so(2n+1)}(2,0))$$ is rational and $C_2$-cofinite if $C(L_{so(2n+1)}(k+1,0),L_{so(2n+1)}(k,0)\otimes L_{so(2n+1)}(1,0))$ and $C(L_{so(2n+1)}(k+2,0),L_{so(2n+1)}(k+1,0)\otimes L_{so(2n+1)}(1,0))$ are rational and $C_2$-cofinite (see Theorem \ref{rationalb}).

  In our study of rationality and $C_2$-cofiniteness of   $C(L_{so(n)}(k+2,0),L_{so(n)}(k,0)\otimes L_{so(n)}(2,0))$ and $C(L_{sl_2}(k+4,0),L_{sl_2}(k,0)\otimes L_{sl_2}(4,0))$, the results of Adamovic and Perse in \cite{AP} play an important role. By using the boson-fermion correspondence, it was shown in \cite{AP} that $C(L_{so(n)}(2,0),L_{so(n)}(1,0)\otimes L_{so(n)}(1,0))$ and $C(L_{sl_2}(4,0),L_{sl_2}(2,0)\otimes L_{sl_2}(2,0))$ are closely related to certain rank-one lattice vertex operator algebras. As a result, we can study rationality and $C_2$-cofiniteness of  $C(L_{so(n)}(2,0),L_{so(n)}(1,0)\otimes L_{so(n)}(1,0))$ and $C(L_{sl_2}(4,0),L_{sl_2}(2,0)\otimes L_{sl_2}(2,0))$ by using the orbifold theory \cite{CM}, \cite{M}.

The paper is organized as follows: In Section 2, we recall some
facts about  vertex superalgebras. In Section 3, we show that the coset vertex operator algebra $C(L_{so(2n)}(k+2,0),L_{so(2n)}(k,0)\otimes L_{so(2n)}(2,0))$ is rational and $C_2$-cofinite when $k$ is an admissible number for $\hat{so(2n)}$. In Section 4, we study rationality and $C_2$-cofiniteness of the coset vertex operator algebra $C(L_{so(2n+1)}(k+2,0),L_{so(2n+1)}(k,0)\otimes L_{so(2n+1)}(2,0))$ under certain assumption. In Section 5, we show that the coset vertex operator algebra $C(L_{sl_2}(k+4,0),L_{sl_2}(k,0)\otimes L_{sl_2}(4,0))$ is rational and $C_2$-cofinite when $k$ is an admissible number for $\hat{sl_2}$. In Section 6, we classify irreducible modules of $C(L_{sl_2}(k+4,0),L_{sl_2}(k,0)\otimes L_{sl_2}(4,0))$ when $k$ is a positive odd integer.

\section{Preliminaries }
\def\theequation{2.\arabic{equation}}
\setcounter{equation}{0}
In this section, we recall some facts about vertex superalgebras. We will continue to use notions in \cite{Lin}. Let $V=V^{even}\oplus V^{odd}$ be a vertex superalgebra (cf. \cite{K2}, \cite{KWang}, \cite{L2}). If $V^{odd}=0$, $V$ is a vertex algebra (cf. \cite{B}). For a subalgebra $U$ of $V$, denote by
 $$C(U, V)=\{v\in V|u_nv=0, \forall u\in U, n\geq 0\}$$ the commutant of $U$ in $V$. Then $C(U, V)$ is a subalgebra of $V$ (also called coset vertex superalgebra) (cf. \cite{FZ}, \cite{LL}).
 \subsection{Tensor product vertex operator algebras} In this subsection, we recall some facts about tensor product vertex operator algebras.  Let $(V^1, Y, \1, \omega^1), \cdots, (V^p, Y, \1,\omega^p)$ be vertex operator algebras. The
tensor product of $(V^1, Y, \1, \omega^1), \cdots, (V^p, Y, \1,\omega^p)$ is
constructed on $V = V^1\otimes
\cdots \otimes V^p$.
The vertex operator  $Y$  of the tensor product is
defined by
$$Y(v^1\otimes \cdots \otimes v^p, z) = Y(v^1, z)\otimes
\cdots \otimes Y(v^p, z)$$for $v^i \in V^i\ (1\leq i \leq p)$.
Then  $(V, Y)$
is a vertex operator algebra such that the vacuum element and Virasoro vector are $ \1\otimes \cdots \otimes \1$
and 
$\omega^1 \otimes \cdots \otimes \1 + \cdots +\1\otimes
\cdots\otimes \omega^p$, respectively  \cite{FHL}, \cite{LL}.

Let $(M^1, Y_{M^1}), \cdots, (M^p, Y_{M^p})$  be  modules of $V^1, \cdots, V^p$, respectively.  We may construct
a $V^1\otimes \cdots \otimes
V^p$-module structure on $M^1\otimes \cdots \otimes M^p$ by
$$Y_{M^1\otimes \cdots \otimes M^p}(v^1\otimes \cdots \otimes v^p,z) =
Y_{M^1}(v^1, z)\otimes \cdots \otimes Y_{M^p}(v^p, z),$$ for $v^i \in V^i\ (1\leq i \leq p).$ Then  $(M^1\otimes \cdots \otimes M^p, Y_{M^1\otimes \cdots \otimes M^p})$ is a $V^1\otimes
\cdots \otimes V^p$-module. The
 following results have been obtained in \cite{FHL}, \cite{DMZ} and \cite{ABD}.
\begin{theorem}\label{tensor}
(i)
$V^1\otimes \cdots \otimes V^p$  is rational if and only if $V^1, \cdots, V^p$ are rational. Moreover, any irreducible $V^1\otimes \cdots \otimes V^p$-module has the form $M^1\otimes \cdots \otimes M^p$, where $M^1, \cdots, M^p$ are
some irreducible modules of $V^1,
\cdots, V^p,$ respectively.\\
 (ii) If $V^1, \cdots, V^p$ are $C_2$-cofinite, then $V^1\otimes \cdots \otimes V^p$  is $C_2$-cofinite.\\
 (iii) If $V^1, \cdots, V^p$ are strongly regular, then $V^1\otimes \cdots \otimes V^p$  is strongly regular.
\end{theorem}

Recall that a vertex operator algebra $V$ is called {\em regular} if every weak $V$-module is a direct sum of irreducible ordinary $V$-modules. Then we have
\begin{proposition}\label{regulart}
Let $V^1$, $V^2$ be simple vertex operator algebras. Suppose that $V^1\otimes V^2$ and $V^2$ are regular, then $V^1$ is regular. In particular, $V^1$ is rational and $C_2$-cofinite.
\end{proposition}
\pf Let $M$ be a weak $V^1$-module. Then $M\otimes V^2$ is a weak $V^1\otimes V^2$-module. Since $V^1\otimes V^2$ and $V^2$ are regular, then $M\otimes V^2$ must be a direct sum of irreducible ordinary $V^1\otimes V^2$-modules $\{ W^s\otimes V^2|s\in S\}$, where $\{W^s|s\in S\}$ are irreducible ordinary $V^1$-modules. This implies $M=\bigoplus_{s\in S} W^s$. Hence, $V^1$ is regular. This implies that $V^1$ is rational and $C_2$-cofinite \cite{Li}.
\qed
\subsection{Fermion vertex operator superalgebras}\label{fermion}

 In this subsection we recall from \cite{K2}, \cite{L2} some facts about fermion vertex operator superalgebras. Let $A$ be a $n$-dimensional vector space with a nondegenerate symmetric bilinear form $\langle,\rangle$. The Clifford affinization of $A$ is a Lie superalgebra
 \begin{align*}
 \hat A=A\otimes \C[t, t^{-1}]t^{\frac{1}{2}}\oplus \C C,
 \end{align*}
 with $\Z_2$-gradation $\hat A_{\bar 0}=\C C$, $\hat A_{\bar 1}=A\otimes \C[t, t^{-1}]t^{\frac{1}{2}}$ and the communication relations:
 \begin{align*}
 [u(m), v(n)]_+=\delta_{m+n, 0}\langle u, v\rangle C,~~[u(m), C]=0,
 \end{align*}
 for any $u, v\in A$, $m, n\in \frac{1}{2}+\Z$, where $u(m)=u\otimes t^m$.

 Set $\hat A_+=t^{\frac{1}{2}}\C[t]\otimes A$ and $\hat A_-=t^{-\frac{1}{2}}\C[t^{-1}]\otimes A$, then $\hat A_+$ and $\hat A_-$ are subalgebras of $\hat A$. For any nonzero complex number $c$, we consider the Verma module
 \begin{align*}
 M_{\hat A}(c, 0)=U(\hat A)/J,
 \end{align*}
  where $U(\hat A)$ denotes the universal enveloping algebra of $\hat{A}$ and $J$ is the left ideal of $U(\hat A)$ generated by $\hat A_+$, $C-c$. Then $M_{\hat A}(c, 0)$ is an irreducible highest weight module for $\hat A$ (cf. \cite{L2}). Moreover,  $M_{\hat A}(c, 0)$ has a vertex operator superalgebra structure (cf. \cite{K2}, \cite{L2}).

\subsection{Lattice vertex operator superalgebras}\label{lattice}
In this subsection we recall from \cite{FLM}, \cite{K2}, \cite{DL} some facts about lattice vertex operator
superalgebras. Let $L$ be a positive definite integral lattice and $(,)$ be the associated positive definite bilinear form.
We
consider the central extension $\hat{L}$ of $L$ by the cyclic
group $\langle \kappa\rangle$ of order $2$:
$$1\to \langle \kappa\rangle\to \hat{L}\to L\to 1,$$ such that the commutator map
$c(\alpha, \beta)=\kappa^{(\alpha, \beta)+(\alpha, \alpha)(\beta, \beta)}$ for any $\alpha, \beta\in L$ \cite{FLM}. Let $e: L\to \hat{L}$
be a section such that $e_0 = 1$ and $\epsilon_0: L\times L\to \langle \kappa\rangle$ be the corresponding 2-cocycle  \cite{FLM}. Then we have
$e_{\alpha}e_{\beta}=\epsilon_0(\alpha, \beta)e_{\alpha+\beta}$ for
$ \alpha, \beta \in L$, and $\epsilon_0(\alpha,
\beta)\epsilon_0(\beta, \alpha)=\kappa^{(\alpha, \beta)+(\alpha, \alpha)(\beta, \beta)}$.

Next, we consider the induced
$\hat{L}$-module:$$\C\{L\}=\C[\hat{L}]\otimes_{\langle \kappa\rangle}\C\cong
\C[L]\ \ (\text{linearly}),$$ where $\C[L]$ denotes the group algebra of $L$ and
$\kappa$ acts on $\C$ as multiplication by $-1$.  Then $\C[L]$ becomes an
$\hat{L}$-module such that $e_{\alpha}\cdot
e^{\beta}=\epsilon(\alpha, \beta)e^{\alpha+\beta}$ and $\kappa\cdot
e^{\beta}=-e^{\beta}$, where $\epsilon(\alpha, \beta)$ is defined as $\nu\circ\epsilon_0$ and $\nu$ is the an isomorphism from $\langle \kappa\rangle$ to $\langle \pm 1\rangle$ such that $\nu(\kappa)=-1$.

 Set $H=\C\otimes_{\Z}L$,  and consider the Heisenberg algebra $\hat{H}=H\otimes \C[t, t^{-1}]\oplus \C K$ with the communication relations: For any $u, v\in H$,
\begin{align*}
&[u(m), v(n)]=m(u, v)\delta_{m+n, 0}K,\\
&[K, x]=0, \text{ for any } x\in \hat{H},
\end{align*}
where $u(n)=u\otimes t^n$, $u\in H$.
Let $M(1)$ be the
Heisenberg vertex operator algebra associated to $\hat H$ (cf. \cite{LL}). The vector space of the lattice vertex operator superalgebra is  defined to be
$$V_L=M(1)\otimes_{\C}\C\{L\}\cong M(1)\otimes_{\C}\C[L]\ \ (\text{linearly}).$$

For an element $h\in H$, we define an action $h(0)$ on
$\C[L]$ by $h(0)\cdot e^{\alpha}=(h, \alpha)e^{\alpha}$ for $\alpha\in {L}$. We also define an action $z^{h(0)}$ on $\C[L]$ by $z^{h(0)}\cdot e^{\alpha}=z^{(h,
\alpha)}e^{\alpha}$. Then $\hat{L}$, $h(n)(n\neq 0)$, $h(0)$ and $z^{h(0)}$ act naturally on $V_L$ by acting on either $M(1)$ or $\C[L]$ as indicated above. It was proved in  \cite{FLM}, \cite{K2} that
$V_L$ has a vertex operator superalgebra structure such that
\begin{align*}
&Y(h(-1)1, z)=h(z)=\sum_{n\in \Z}h(n)z^{-n-1}\ \ (h\in H),\\
&Y(e^{\alpha}, z)=E^-(-\alpha, z)E^+(-\alpha, z)e_{\alpha}z^{\alpha(0)},
\end{align*}
where
\begin{align*}
E^-(\alpha, z)=\exp\left(\sum_{n<0}\frac{\alpha(n)}{n}z^{-n}\right),\ \ \
E^+(\alpha, z)=\exp\left(\sum_{n>0}\frac{\alpha(n)}{n}z^{-n}\right).
\end{align*}

Denote by $L^\circ$ the dual lattice of $L$. Then $V_{L+\gamma}$, for $\gamma\in L^\circ/L$ are all irreducible modules of $V_L$ (cf. \cite{DL}).
\subsection{Affine vertex operator algebras}
In this subsection, we shall recall some facts about affine vertex operator algebras from \cite{K}, \cite{LL}.
Let $\g$ be a finite dimensional simple Lie algebra and $(\, ,\, )$ the normalized Killing form of $\g$, i.e., $(\theta, \theta)=2$ for the highest root $\theta$ of $\g$. Fix a Cartan subalgebra $\h$ of $\g$ and  denote the corresponding root system by $\Delta_{\g}$ and the root lattice by $Q$. We further fix simple roots $\{\alpha_1,\cdots,\alpha_n\}\subset\h^*$ and simple coroots $\{\alpha_1^\vee,\cdots,\alpha_n^\vee\}\subset\h$. Then the weight lattice $P$ of $\g$ is the set of $\lambda\in \h^*$ such that $\langle\lambda, \alpha_i^\vee\rangle\in\Z$ for  $i=1,\cdots, n$. Note that $P$ is equal to $\oplus_{i=1}^n\Z\Lambda_i$, where $\Lambda_i$ are the fundamental weights defined by the equation $\langle\lambda_i, \alpha_j^\vee\rangle=\delta_{i,j}$ for  $i,j=1,\cdots, n$. We also use the standard notation $P_+$ to denote the set of dominant weights $\{\Lambda\in P\mid\langle\lambda, \alpha_i^\vee\rangle\geq 0,~1\leq i\leq n \}$.

Recall that the affine Lie algebra of $\g$ is defined on $\hat{\g}=\g\otimes \C[t^{-1}, t]\oplus \C K$ with Lie brackets
\begin{align*}
[x(m), y(n)]&=[x, y](m+n)+(x, y)m\delta_{m+n,0}K,\\
[K, \hat\g]&=0,
\end{align*}
for $x, y\in \g$ and $m,n \in \Z$, where $x(n)$ denotes $x\otimes t^n$.

For a complex number $k$ and $\Lambda \in P$, let $L_{\g}(\Lambda)$ be the irreducible highest weight module for $\g$ with highest weight $\Lambda$ and define
\begin{align*}
V_{\g}(k, \Lambda)=\Ind_{\g\otimes \C[t]\oplus \C K}^{\hat \g}L_{\g}(\Lambda),
\end{align*}
where $L_{\g}(\Lambda)$ is viewed as a module for $\g\otimes \C[t]\oplus \C K$ such that $\g\otimes t\C[t]$ acts as $0$ and $K$ acts as $k$. It is well-known that $V_{\g}(k, \Lambda)$ has a unique maximal proper submodule which is denoted by $J(k, \Lambda)$ (see \cite{K}). Let $L_{\g}(k, \Lambda)$ be the corresponding irreducible quotient module.  It was proved in \cite{FZ} that $L_{\g}(k, 0)$ has a vertex operator algebra structure if $k\neq -h^\vee$, where $h^\vee$ denotes the dual Coxeter number of $\g$. Moreover, the following results have been obtained in \cite{DLM2}, \cite{FZ}.
\begin{theorem}\label{moduleaff}
Let $k$ be a positive integer. Then\\
 (1) $L_{\g}(k, 0)$ is a strongly regular vertex operator algebra.\\
 (2) $L_{\g}(k, \Lambda)$ is a module for the vertex operator algebra $L_{\g}(k, 0)$ if and only if $\Lambda \in P_+^k$, where $P_+^k=\{\Lambda \in P_+|(\Lambda, \theta)\leq k\}$.
\end{theorem}

\subsection{Diagonal coset vertex operator algebras}
 In this subsection, we assume that $\g$ is simply laced. Recall that the level $k$ of $L_{\g}(k,0)$ is called an {\em admissible number} for $\hat \g$ if  $L_{\g}(k,0)$ is an admissible module of $\hat \g$ \cite{KW1}. Since  $\g$ is simply laced, this condition is equivalent to that
 $$k+h^\vee=\frac{p}{q},~p, q\in \Z_{\geq 1},~(p, q)=1,~p\geq h^\vee.$$

 For an admissible number $k$ for $\hat \g$  and a positive integer $l$, we consider the vertex operator algebra $L_{\g}(k,0)\otimes L_{\g}(l,0)$. Let $W$ be the vertex subalgebra of  $L_{\g}(k,0)\otimes L_{\g}(l,0)$ generated by $\{x(-1)\1\otimes \1+\1\otimes x(-1)\1|x\in \g\}$. By Corollary 4.1 of \cite{KW1}, $W$ is isomorphic to $L_{\g}(k+l,0)$. Set \begin{align*}
C(L_{\g}&(k+l,0),L_{\g}(k,0)\otimes L_{\g}(l,0))
\\&=\{u\in L_{\g}(k,0)\otimes L_{\g}(l,0)|u_{n}v=0, \forall v\in L_{\g}(k+l,0), \forall n\in \Z_{\geq 0}\}.
\end{align*}
It is well-known that $C(L_{\g}(k+l,0),L_{\g}(k,0)\otimes L_{\g}(l,0))$ is a vertex subalgebra of $L_{\g}(k,0)\otimes L_{\g}(l,0)$ \cite{LL}. Moreover, we have the following result.
\begin{proposition}\label{scoset}
Let $\w^1, \w^2, \w^a$ be the Virasoro vectors of $L_{\g}(k,0)$, $L_{\g}(l,0)$, $L_{\g}(k+l,0)$, respectively. Then $C(L_{\g}(k+l,0),L_{\g}(k,0)\otimes L_{\g}(l,0))$  is a  simple vertex operator algebra with the Virasoro vector $\w^1+\w^2-\w^a$.
\end{proposition}
\pf By the discussion above, we know that $W$ is a simple vertex operator subalgebra of $L_{\g}(k,0)\otimes L_{\g}(l,0)$.  Moreover, by Corollary 4.1 of \cite{KW1}, $L_{\g}(k,0)\otimes L_{\g}(l,0)$ viewed as a module of $L_{\g}(k+l,0)$ is completely reducible. Then it follows from Lemma 2.1 of \cite{ACKL} that $C(L_{\g}(k+l,0),L_{\g}(k,0)\otimes L_{\g}(l,0))$  is a  simple vertex  algebra. Furthermore, it follows from Theorem 3.11.12 of \cite{LL} that $C(L_{\g}(k+l,0),L_{\g}(k,0)\otimes L_{\g}(l,0))$  is a  simple vertex operator algebra with the Virasoro vector $\w^1+\w^2-\w^a$.
\qed

When $l=1$, the following  result about rationality and $C_2$-cofiniteness of $C(L_{\g}(k+1,0),L_{\g}(k,0)\otimes L_{\g}(1,0))$ has been established in \cite{ACL}.
\begin{theorem}\label{rational}
Let $\g$ be a simply laced simple Lie algebra, $k$ be an admissible number for $\hat \g$. Then the vertex operator algebra $C(L_{\g}(k+1,0),L_{\g}(k,0)\otimes L_{\g}(1,0))$ is strongly regular.
\end{theorem}
\begin{remark}
When $k=1$, rationality and $C_2$-cofiniteness of  the vertex operator algebra $C(L_{\g}(k+1,0),L_{\g}(k,0)\otimes L_{\g}(1,0))$ have also been established for the case $\g=sl_{n}$ in \cite{ACL}, \cite{JL1}, \cite{JL2} and the case $\g=so(n)$ in \cite{AP}.
\end{remark}

\section{Rationality and $C_2$-cofiniteness of diagonal coset vertex operator algebras $C(L_{\g}(k+2,0),L_{\g}(k,0)\otimes L_{\g}(2,0))$: type $D_n$}
\def\theequation{3.\arabic{equation}}
\setcounter{equation}{0}
In this section, we shall prove that  the coset vertex operator algebra $C(L_{so(2n)}(k+2,0),L_{so(2n)}(k,0)\otimes L_{so(2n)}(2,0))$ is rational and $C_2$-cofinite if $n\geq 3$ and  $k$ is an admissible number for $\hat{so(2n)}$.
\subsection{The commutant of   $L_{so(2n)}(2,0)$ in $L_{so(2n)}(1,0)\otimes L_{so(2n)}(1,0)\oplus L_{so(2n)}(1,\Lambda_1)\otimes L_{so(2n)}(1,\Lambda_1)$} In this subsection, we will show that the commutant of   $L_{so(2n)}(2,0)$ in the vertex operator algebra $L_{so(2n)}(1,0)\otimes L_{so(2n)}(1,0)\oplus L_{so(2n)}(1,\Lambda_1)\otimes L_{so(2n)}(1,\Lambda_1)$ is isomorphic to a rank-one lattice vertex operator algebra. This fact will play an important role in  our study of rationality and $C_2$-cofiniteness of $C(L_{so(2n)}(k+2,0),L_{so(2n)}(k,0)\otimes L_{so(2n)}(2,0))$.

We first recall  some results in \cite{AP}. Let $A$ be a vector space spanned by $\Psi^{\pm}_{i}$,  $\Phi^{\pm}_{i}$, $i=1, \cdots, n$, and $(,)$ be the nondegenerate symmetric bilinear form on $A$ such that
$$(\Psi^{\pm}_i, \Psi^{\mp}_j)=(\Phi^{\pm}_i, \Phi^{\mp}_j)=\delta_{i,j}, ~~~(\Psi^{\pm}_i, \Psi^{\pm}_j)=(\Phi^{\pm}_i, \Phi^{\pm}_j)=0,$$
$$(\Psi^{\pm}_i, \Phi^{\pm}_j)=(\Psi^{\mp}_i, \Phi^{\pm}_j)=0,$$
for $i,j=1, \cdots, n$. Then we have a fermion vertex operator superalgebra $M_{\hat A}(1, 0)$ associated to $\hat A$ (see Subsection \ref{fermion}). In the following, we will use $F_{2n}$ to denote the fermion vertex operator superalgebra $M_{\hat A}(1, 0)$.

Let $F^{\Psi}_n$ (resp. $F^{\Phi}_n$) be the subalgebra of $F_{2n}$ generated by $\Psi_i^{\pm}(-\frac{1}{2})\1$ (resp. $\Phi_i^{\pm}(-\frac{1}{2})\1$), $i=1, \cdots, n$. Then the following results are well-known (see \cite{Fe}, \cite{FF}).
\begin{theorem}\label{FD}
Let $n$ be a positive integer such that $n\geq 3$. Then\\
 (i) The vertex operator algebras $(F_n^{\Psi})^{even}$ and $(F_n^{\Phi})^{even}$ are isomorphic to $L_{so(2n)}(1,0)$.
 (ii) $(F_n^{\Psi})^{odd}$ and $(F_n^{\Phi})^{odd}$ viewed as $L_{so(2n)}(1,0)$-modules are isomorphic to $L_{so(2n)}(1,\Lambda_1)$.
\end{theorem}

Define the  lattice
$$R_{2n}=\Z x_1+\cdots+\Z x_n+\Z y_1+\cdots+\Z y_n,$$
$$(x_i, x_j)=(y_i, y_j)=\delta_{i,j},~~~(x_i, y_j)=0,$$
where $i, j=1, \cdots, n$. We set $z_{2k-1}=x_k, z_{2k}=y_k$, where $k=1, \cdots, n$.

 Following \cite{AP},  we choose the 2-cocycle $\epsilon: R_{2n}\times R_{2n}\to \{\pm1\}$ determined by
 \begin{align}
 \epsilon(z_i, z_j)=\left\{
\begin{array}{ll}
1,& \text{~if~}i\leq j,\\
-1,& \text{~if~} i>j.
\end{array}
\right.
 \end{align}
 Then we  have a lattice vertex operator superalgebra $V_{R_{2n}}$ associated to $R_{2n}$ (see Subsection \ref{lattice}). We shall need the following (non-standard) version of the boson-fermion correspondence \cite{Fe}, \cite{F}, \cite{K2}  (see also Theorem 2 of \cite{AP}).
 \begin{theorem}\label{bf}
 Let $n$ be a positive integer such that $n\geq 3$. Then there exists a vertex superalgebra isomorphism $\phi_{2n}: F_{2n}\to V_{R_{2n}}$ such that
 \begin{align*}
 &\Psi^+_k(-\frac{1}{2})\1\mapsto \frac{1}{\sqrt{2}}(e^{x_k}+e^{-y_k}),\ \ \
 \Psi^-_k(-\frac{1}{2})\1\mapsto \frac{1}{\sqrt{2}}(e^{-x_k}+e^{y_k}),\\
 &\Phi^+_k(-\frac{1}{2})\1\mapsto \frac{i}{\sqrt{2}}(e^{x_k}-e^{-y_k}),\ \ \
 \Phi^-_k(-\frac{1}{2})\1\mapsto \frac{-i}{\sqrt{2}}(e^{-x_k}-e^{y_k}),
 \end{align*}
 for $k=1,\cdots, n$.
 \end{theorem}

 Based on Theorems \ref{FD} and \ref{bf}, the following results have been obtained in Theorem 3 of \cite{AP}.
 \begin{theorem}\label{subD}
 Let $n$ be a positive integer such that $n\geq 3$. Then\\
 (i) The subalgebra of $V_{R_{2n}}$ generated by
 \begin{align}\label{generator1}
 e^\gamma \text{~ and~} e^{-\gamma},
  \end{align}where $\gamma=x_1+\cdots+x_n+y_1+\cdots+y_n$, is isomorphic to the rank-one lattice vertex operator algebra $V_L$ such that $$L=\Z\gamma {\text~ and~}(\gamma, \gamma)=2n.$$
 (ii)  The subalgebra of $V_{R_{2n}}$ generated by
 \begin{align}\label{generator2}
 e^{x_k-y_l}+e^{x_l-y_k},~~~e^{x_k-x_l}+e^{y_l-y_k}, ~~~e^{-x_k+y_l}+e^{-x_l+y_k},
 \end{align}
 for $k, l=1, \cdots, n$ and $k\neq l$, is isomorphic to $L_{so(2n)}(2, 0)$.\\
 (iii) The vertex operator algebra $V_L\otimes L_{so(2n)}(2, 0)$ is isomorphic to the subalgebra of $V_{R_{2n}}$ generated by elements (\ref{generator1}) and (\ref{generator2}).
 \end{theorem}

 Furthermore, it was noticed in \cite{AP} that the root lattice of the Lie algebra $sl_{2n}$ can be realized as a sublattice of $R_{2n}$:
 $$A_{2n-1}=\Z(x_1-x_2)+\cdots+\Z(x_{n-1}-x_n)+\Z(x_n-y_n)+\Z(y_n-y_{n-1})+\cdots+\Z(y_2-y_1).$$
Let $\lambda_0=0$ and
\begin{align*}
&\lambda_i=x_1+\cdots+ x_i,~~~\text{~for~} 1\leq i\leq n,\\
&\lambda_{n+i}=x_1+\cdots+ x_n+y_n+\cdots+ y_{n-i+1},~~~\text{~for~} 1\leq i\leq n-1.
\end{align*}
Then the following results have been obtained in \cite{AP}.
\begin{theorem}\label{decomp1}
Let $n$ be a positive integer such that $n\geq 3$. Then\\
(i) $A_{2n-1}=L^{\perp}$.\\
(ii) $V_{R_{2n}}$ viewed as a module of $V_L\otimes L_{sl_{2n}}(1, 0)$ has the following decomposition
\begin{align}
V_{R_{2n}}\cong\bigoplus_{i=0}^{2n-1}V_{\lambda_i+L+A_{2n-1}}\cong \bigoplus_{i=0}^{2n-1} V_{L+\frac{i}{2n}\gamma}\otimes L_{sl_{2n}}(1, \Lambda_i).
\end{align}
\end{theorem}

By Theorems \ref{subD}, \ref{decomp1}, the affine vertex operator algebra $L_{so(2n)}(2, 0)$ is a subalgebra of $L_{sl_{2n}}(1, 0)$ (see Remark 1 of \cite{AP}). Moreover, the following result has been obtained in \cite{KW}  (see also Example 7.3.4 of \cite{W}).
\begin{proposition}\label{decompwaki}
Let $n$ be a positive integer such that $n\geq 3$. Then $L_{sl_{2n}}(1, 0)$ viewed as a module of $L_{so(2n)}(2, 0)$ has the following decomposition
$$L_{sl_{2n}}(1, 0)\cong L_{so(2n)}(2, 0)\oplus L_{so(2n)}(2, 2\Lambda_1).$$
\end{proposition}

Combining Theorems \ref{FD}, \ref{decomp1}, we have
\begin{proposition}\label{even}
Let $n$ be a positive integer such that $n\geq 3$. Then\\
(i) The even part $F_{2n}^{even}$ of $F_{2n}$ contains a subalgebra isomorphic to $L_{so(2n)}(1, 0)\otimes L_{so(2n)}(1, 0)$. Moreover, $F_{2n}^{even}$ viewed as a module of  $L_{so(2n)}(1, 0)\otimes L_{so(2n)}(1, 0)$ is isomorphic to $L_{so(2n)}(1, 0)\otimes L_{so(2n)}(1, 0)\oplus L_{so(2n)}(1, \Lambda_1)\otimes L_{so(2n)}(1, \Lambda_1)$.\\
(ii)  The even part $V_{R_{2n}}^{even}$ of $V_{R_{2n}}$ contains the subalgebra $V_L\otimes L_{sl_{2n}}(1, 0)$. Moreover, $V_{R_{2n}}^{even}$ viewed as a module of $V_L\otimes L_{sl_{2n}}(1, 0)$ has the following decomposition
\begin{align*}
V_{R_{2n}}^{even}\cong \bigoplus_{\substack{0\leq i\leq n;\\i: \text{even}}}V_{L+\frac{i}{2n}\gamma}\otimes L_{sl_{2n}}(1, \Lambda_i)\bigoplus \bigoplus_{\substack{1\leq i\leq n-1;\\n+i: \text{even}}}V_{L+\frac{n+i}{2n}\gamma}\otimes L_{sl_{2n}}(1, \Lambda_{n+i}).
\end{align*}
\end{proposition}

We are now ready to determine the commutant of   $L_{so(2n)}(2,0)$ in the vertex operator algebra $L_{so(2n)}(1,0)\otimes L_{so(2n)}(1,0)\oplus L_{so(2n)}(1,\Lambda_1)\otimes L_{so(2n)}(1,\Lambda_1)$.
\begin{theorem}\label{cosetD}
Let $n$ be a positive integer such that $n\geq 3$.
Then the commutant of   $L_{so(2n)}(2,0)$ in  $L_{so(2n)}(1,0)\otimes L_{so(2n)}(1,0)\oplus L_{so(2n)}(1,\Lambda_1)\otimes L_{so(2n)}(1,\Lambda_1)$ is isomorphic to the lattice vertex algebra $V_L$.
\end{theorem}
\pf  By Theorem \ref{bf} and Propositions \ref{decompwaki}, \ref{even}, we have
\begin{align*}
&C(L_{so(2n)}(2,0), L_{so(2n)}(1,0)\otimes L_{so(2n)}(1,0)\oplus L_{so(2n)}(1,\Lambda_1)\otimes L_{so(2n)}(1,\Lambda_1))\\
&\cong C(L_{so(2n)}(2,0), F_{2n}^{even})\\
&\cong C(L_{so(2n)}(2,0), R_{2n}^{even})\cong V_L,
\end{align*}
as desired.
\qed
\subsection{Rationality and $C_2$-cofiniteness of diagonal coset vertex operator algebras $C(L_{so(2n)}(k+2,0),L_{so(2n)}(k,0)\otimes L_{so(2n)}(2,0))$} We first recall some facts about orbifold vertex operator algebras. Let $(V, Y, 1, \omega)$ be a vertex operator algebra. An {\em automorphism} of $V$ is a linear isomorphism $\sigma: V\to V$ such that:
$$\sigma(u_nv)=\sigma(u)_n\sigma(v),~~~\sigma(\w)=\w,~~~\sigma(\1)=\1,$$ for any $n\in \Z$ and $u, v\in V$.  For a finite automorphism group $G$ of $V$, we define $V^G=\{v\in V|g(v)=v, \forall g\in G\}$. Then $V^G$ is a subalgebra of $V$. We will need the following result, which has been obtained in \cite{CM}, \cite{M}.
\begin{theorem}\label{orbifold}
Let $V$ be a $C_2$-cofinite simple vertex operator algebra of CFT type, $G$ be a finite automorphism group of $V$. Suppose that $G$ is abelian, then the vertex operator algebra  $V^G$ is $C_2$-cofinite. Furthermore, $V^G$ is rational if $V$ is strongly regular and $G$ is abelian.
\end{theorem}

As an application of Theorem \ref{orbifold}, we have
\begin{lemma}\label{cosetr}
Let $V$ be a strongly regular vertex operator algebra. Suppose that \\
(i) $V$ contains a rank-one lattice vertex operator algebra  $V_{L}$ with $L=\Z\beta$.\\
(ii) $C(V_L, V)$ is a vertex operator algebra such that the Virasoro vector is $\w-\w_L$, where $\w$ and $\w_L$ denote the Virasoro vectors of $V$ and $V_L$, respectively.\\
 Then the vertex operator algebra $C(V_L, V)\otimes V_L$ is rational and $C_2$-cofinite.
\end{lemma}
\pf Denote by $(,)$ the bilinear form of $L$. We then assume that $(\beta, \beta)=m$ for some positive even integer $m$.  Since $V_L$ is a regular vertex operator algebra \cite{DLM1}, $V$ viewed as a $V_L$-module has the following decomposition
$$V=\bigoplus_{0\leq r\leq m-1}M^r\otimes V_{L+\frac{r}{m}\beta},$$
where $M^r$ denotes the multiplicity space of $V_{L+\frac{r}{m}\beta}$ in $V$. Since each irreducible $V_L$-module $V_{L+\frac{r}{m}\beta}$ is a simple current $V_L$-module \cite{DLM0}, the set $\mathbf{G}=\{r|0\leq r\leq m-1,~M^r\neq 0\}$ forms a subgroup of $\Z/\Z_{m}$. Let $s$ be the smallest nonzero integer such that $M^s\neq 0$. Then $s$ is a generator of $\mathbf G$ and we have
$$V=\bigoplus_{j:~0\leq js\leq m-1}M^{js}\otimes V_{L+\frac{js}{m}\beta}.$$
Denote by $t$ the order of $\mathbf G$ and  set $\zeta=e^{2\pi i /t}$. We then define a linear map $g$ of $V$ such that $g|_{M^{js}\otimes V_{L+\frac{js}{m}\beta}}=\zeta^j\id$. Since each irreducible $V_L$-module $V_{L+\frac{r}{m}\beta}$ is a simple current $V_L$-module, $g$ is an automorphism of the vertex operator algebra $V$. Moreover, $V^{\langle g\rangle}=M^0\otimes V_L=C(V_L, V)\otimes V_L$. Hence, by Theorem \ref{orbifold},  $C(V_L, V)\otimes V_L$ is rational and $C_2$-cofinite.
\qed

\vskip.25cm
  We next recall some facts about extensions of vertex operator algebras. Recall that a vertex operator algebra $U$ is called an {\em extension} of $V$ if $V$ is a vertex subalgebra of $U$, and $V$, $U$ have the same Virasoro vector. We will  need the following result about extensions which has been obtained in \cite{ABD}.
    \begin{theorem}\label{abd}
Let $V$ be a $C_2$-cofinite vertex operator algebra and $U$ be an extension of $V$. Assume further that $U$ viewed as a $V$-module is completely reducible, then $U$ is $C_2$-cofinite.
\end{theorem}

We also  need the following result about extensions which has been obtained in  \cite{HKL}.
  \begin{theorem}\label{extensionr}
Let $V$ be a strongly regular vertex operator algebra. Suppose that $U$ is a simple vertex operator algebra and an extension of $V$, then $U$ is rational.
\end{theorem}

We are now ready to prove  the main result in this section.
\begin{theorem}\label{rationalD}
 Let $n$ be a positive integer such that $n\geq 3$ and $k$ be an admissible number for $\hat{so(2n)}$. Then the  vertex operator algebra $C(L_{so(2n)}(k+2,0),L_{so(2n)}(k,0)\otimes L_{so(2n)}(2,0))$ is rational and $C_2$-cofinite.
\end{theorem}
\pf Consider the vertex operator algebra $L_{so(2n)}(k,0)\otimes F_{2n}^{even}$. By Proposition \ref{even}, $L_{so(2n)}(k,0)\otimes F_{2n}^{even}$ contains a subalgebra $L_{so(2n)}(k,0)\otimes L_{so(2n)}(1,0)\otimes L_{so(2n)}(1,0)$. Moreover,  $L_{so(2n)}(k,0)\otimes F_{2n}^{even}$ viewed as a module of $L_{so(2n)}(k,0)\otimes L_{so(2n)}(1,0)\otimes L_{so(2n)}(1,0)$ has the following decomposition $$ L_{so(2n)}(k,0)\otimes L_{so(2n)}(1,0)\otimes L_{so(2n)}(1,0)\oplus L_{so(2n)}(k,0)\otimes L_{so(2n)}(1,\Lambda_1)\otimes L_{so(2n)}(1,\Lambda_1).$$
Therefore, by Theorem \ref{rational}, $L_{so(2n)}(k,0)\otimes F_{2n}^{even}$ is an extension of
\begin{align*}
L_{so(2n)}&(k+2,0)\otimes C(L_{so(2n)}(k+2,0), L_{so(2n)}(k+1,0)\otimes L_{so(2n)}(1,0))\\
&\otimes C(L_{so(2n)}(k+1,0), L_{so(2n)}(k,0)\otimes L_{so(2n)}(1,0)).
 \end{align*}
 As a result,  $C(L_{so(2n)}(k+2,0), L_{so(2n)}(k,0)\otimes F_{2n}^{even})$ is an extension of
\begin{align*}
C(L_{so(2n)}(k+2,0), &L_{so(2n)}(k+1,0)\otimes L_{so(2n)}(1,0))\\
&\otimes C(L_{so(2n)}(k+1,0), L_{so(2n)}(k,0)\otimes L_{so(2n)}(1,0)).
 \end{align*}
It follows from  Theorems \ref{tensor}, \ref{rational}, \ref{abd} and \ref{extensionr} that $$C(L_{so(2n)}(k+2,0), L_{so(2n)}(k,0)\otimes F_{2n}^{even})$$ is rational and $C_2$-cofinite. By Theorem \ref{bf}, $F_{2n}^{even}$ is isomorphic to $V_{R_{2n}}^{even}$. Hence, $C(L_{so(2n)}(k+2,0), L_{so(2n)}(k,0)\otimes V_{R_{2n}}^{even})$ is rational and $C_2$-cofinite.

We next show that $C(L_{so(2n)}(k+2,0), L_{so(2n)}(k,0)\otimes V_{R_{2n}}^{even})$ is strongly regular. Since $L_{so(2n)}(k,0)\otimes V_{R_{2n}}^{even}$ is of CFT type, $C(L_{so(2n)}(k+2,0), L_{so(2n)}(k,0)\otimes V_{R_{2n}}^{even})$ is of CFT type. Note that $L_{so(2n)}(k,0)\otimes V_{R_{2n}}^{even}$ contains a simple vertex subalgebra isomorphic to $L_{so(2n)}(k+2,0)$. Moreover, $L_{so(2n)}(k,0)\otimes V_{R_{2n}}^{even}$ viewed as an $L_{so(2n)}(k+2,0)$-module is completely reducible. Then it follows from Lemma 2.1 of \cite{ACKL} that $$C(L_{so(2n)}(k+2,0), L_{so(2n)}(k,0)\otimes V_{R_{2n}}^{even})$$ is a simple vertex operator algebra. By Corollary 3.2 of \cite{L1}, this implies that $$C(L_{so(2n)}(k+2,0), L_{so(2n)}(k,0)\otimes V_{R_{2n}}^{even})$$ is self-dual. Thus, $C(L_{so(2n)}(k+2,0), L_{so(2n)}(k,0)\otimes V_{R_{2n}}^{even})$ is strongly regular.

On the other hand, by Propositions \ref{decompwaki}, \ref{even},  $L_{so(2n)}(k,0)\otimes V_{R_{2n}}^{even}$ contains a subalgebra isomorphic to $L_{so(2n)}(k,0)\otimes L_{so(2n)}(2,0)\otimes V_{L}$. Hence, $L_{so(2n)}(k,0)\otimes V_{R_{2n}}^{even}$ is an extension of  $$L_{so(2n)}(k+2,0)\otimes C(L_{so(2n)}(k+2,0), L_{so(2n)}(k,0)\otimes L_{so(2n)}(2,0))\otimes V_{L}.$$ Since $C(L_{so(2n)}(k+2,0), L_{so(2n)}(k,0)\otimes V_{R_{2n}}^{even})$ is strongly regular,  by Lemma \ref{cosetr}, we have
\begin{align*}
&C(V_L, C(L_{so(2n)}(k+2,0), L_{so(2n)}(k,0)\otimes V_{R_{2n}}^{even}))\otimes V_L\\
&=C(L_{so(2n)}(k+2,0)\otimes V_L, L_{so(2n)}(k,0)\otimes V_{R_{2n}}^{even})\otimes V_L
\end{align*}
is rational and $C_2$-cofinite. Moreover, $C(L_{so(2n)}(k+2,0)\otimes V_L, L_{so(2n)}(k,0)\otimes V_{R_{2n}}^{even})\otimes V_L$  is of CFT type, it follows from Theorem 4.5 of \cite{ABD} that $$C(L_{so(2n)}(k+2,0)\otimes V_L, L_{so(2n)}(k,0)\otimes V_{R_{2n}}^{even})\otimes V_L$$  is regular. It is known \cite{DLM1} that $V_L$ is a regular vertex operator algebra. Hence,  by Proposition \ref{regulart}, $C(L_{so(2n)}(k+2,0)\otimes V_L, L_{so(2n)}(k,0)\otimes V_{R_{2n}}^{even})$ is regular. In particular, $C(L_{so(2n)}(k+2,0)\otimes V_L, L_{so(2n)}(k,0)\otimes V_{R_{2n}}^{even})$ is rational and $C_2$-cofinite.

We next show that $C(L_{so(2n)}(k+2,0)\otimes V_L, L_{so(2n)}(k,0)\otimes V_{R_{2n}}^{even})$ is strongly regular. Since $ L_{so(2n)}(k,0)\otimes V_{R_{2n}}^{even}$ is of CFT type, $C(L_{so(2n)}(k+2,0)\otimes V_L, L_{so(2n)}(k,0)\otimes V_{R_{2n}}^{even})$ is  of CFT type. Note that $L_{so(2n)}(k,0)\otimes V_{R_{2n}}^{even}$ contains a simple vertex subalgebra  isomorphic to $L_{so(2n)}(k+2,0)\otimes V_L$. Moreover, $L_{so(2n)}(k,0)\otimes V_{R_{2n}}^{even}$ viewed as a module of $L_{so(2n)}(k+2,0)\otimes V_L$ is completely reducible. Then it follows from Lemma 2.1 of \cite{ACKL} that $C(L_{so(2n)}(k+2,0)\otimes V_L, L_{so(2n)}(k,0)\otimes V_{R_{2n}}^{even})$ is a simple vertex operator algebra. By Corollary 3.2 of \cite{L1}, this implies that $C(L_{so(2n)}(k+2,0)\otimes V_L, L_{so(2n)}(k,0)\otimes V_{R_{2n}}^{even})$ is self-dual. As a consequence, $C(L_{so(2n)}(k+2,0)\otimes V_L, L_{so(2n)}(k,0)\otimes V_{R_{2n}}^{even})$ is strongly regular.

Note that $C(L_{so(2n)}(k+2,0)\otimes V_L, L_{so(2n)}(k,0)\otimes V_{R_{2n}}^{even})$ is an extension of $$C(L_{so(2n)}(k+2,0), L_{so(2n)}(k,0)\otimes L_{so(2n)}(2,0)).$$ We now determine the decomposition of  $C(L_{so(2n)}(k+2,0)\otimes V_L, L_{so(2n)}(k,0)\otimes V_{R_{2n}}^{even})$ as a $C(L_{so(2n)}(k+2,0), L_{so(2n)}(k,0)\otimes L_{so(2n)}(2,0))$-module. By Propositions \ref{decompwaki}, \ref{even}, $C(L_{so(2n)}(k+2,0)\otimes V_L, L_{so(2n)}(k,0)\otimes V_{R_{2n}}^{even})$ has the following decomposition
\begin{align}\label{decompcosetD}
&C(L_{so(2n)}(k+2,0)\otimes V_L, L_{so(2n)}(k,0)\otimes V_{R_{2n}}^{even})\notag\\
&=C(L_{so(2n)}(k+2,0), L_{so(2n)}(k,0)\otimes L_{so(2n)}(2,0))\notag\\
&\ \ \ \oplus C(L_{so(2n)}(k+2,0), L_{so(2n)}(k,0)\otimes L_{so(2n)}(2,2\Lambda_1)).
\end{align}

To show that $C(L_{so(2n)}(k+2,0), L_{so(2n)}(k,0)\otimes L_{so(2n)}(2,0))$ is rational and $C_2$-cofinite, we define a linear map $\tilde\sigma\in \End(L_{so(2n)}(k,0)\otimes L_{sl_{2n}}(1,0))$ such that
$$\tilde\sigma|_{L_{so(2n)}(k,0)\otimes L_{so(2n)}(2,0)}=\id, \ \ \  \tilde\sigma|_{L_{so(2n)}(k,0)\otimes L_{so(2n)}(2,2\Lambda_1)}=-\id.$$
Note that $L_{so(2n)}(2,2\Lambda_1)$ is a simple current module of $L_{so(2n)}(2,0)$ (see Remark 2 of \cite{AP}). Therefore, $\tilde\sigma$ is an automorphism of the vertex operator algebra  $L_{so(2n)}(k,0)\otimes L_{sl_{2n}}(1,0)$. Denote by $\sigma$ the restriction of $\tilde\sigma$ on \begin{align*}
&C(L_{so(2n)}(k+2,0), L_{so(2n)}(k,0)\otimes L_{sl_{2n}}(1,0))\\
&=C(L_{so(2n)}(k+2,0), L_{so(2n)}(k,0)\otimes L_{so(2n)}(2,0)\oplus L_{so(2n)}(k,0)\otimes L_{so(2n)}(2,2\Lambda_1)).
\end{align*}
Then, by the formula (\ref{decompcosetD}), $\sigma$ is an automorphism of $$C(L_{so(2n)}(k+2,0)\otimes V_L, L_{so(2n)}(k,0)\otimes V_{R_{2n}}^{even})$$ such that
\begin{align*}
C(L_{so(2n)}(k+2,0)\otimes V_L, &L_{so(2n)}(k,0)\otimes V_{R_{2n}}^{even})^{\langle\sigma\rangle}\\
&=C(L_{so(2n)}(k+2,0), L_{so(2n)}(k,0)\otimes L_{so(2n)}(2,0)).
\end{align*}
We have proved that $C(L_{so(2n)}(k+2,0)\otimes V_L, L_{so(2n)}(k,0)\otimes V_{R_{2n}}^{even})$ is strongly regular. Then it follows from Theorem \ref{orbifold} that $C(L_{so(2n)}(k+2,0), L_{so(2n)}(k,0)\otimes L_{so(2n)}(2,0))$ is rational and $C_2$-cofinite.
\qed

\vskip.25cm
Since the simple Lie algebra $so(6)$ is isomorphic to $sl_4$, by Theorem \ref{rationalD}, we have
\begin{corollary}\label{rationala3}
 Let $k$ be an admissible number for $\hat{sl_4}$. Then the  vertex operator algebra $C(L_{sl_4}(k+2,0),L_{sl_4}(k,0)\otimes L_{sl_4}(2,0))$ is rational and $C_2$-cofinite.
\end{corollary}

\section{Diagonal coset vertex operator algebras $C(L_{\g}(k+2,0),L_{\g}(k,0)\otimes L_{\g}(2,0))$: type $B_n$}
\def\theequation{4.\arabic{equation}}
\setcounter{equation}{0}
In this section, we shall study rationality and $C_2$-cofiniteness of  diagonal coset vertex operator algebras  $C(L_{so(2n+1)}(k+2,0), L_{so(2n+1)}(k,0)\otimes L_{so(2n+1)}(2,0))$ under certain assumption.
\subsection{The commutant of   $L_{so(2n+1)}(2,0)$ in  vertex operator algebra $L_{so(2n+1)}(1,0)\otimes L_{so(2n+1)}(1,0)\oplus L_{so(2n+1)}(1,\Lambda_1)\otimes L_{so(2n+1)}(1,\Lambda_1)$} In this subsection, we will show that the commutant of   $L_{so(2n+1)}(2,0)$ in $L_{so(2n+1)}(1,0)\otimes L_{so(2n+1)}(1,0)\oplus L_{so(2n+1)}(1,\Lambda_1)\otimes L_{so(2n+1)}(1,\Lambda_1)$ is isomorphic to a rank-one lattice vertex operator algebra. This fact will play an important role in  our study of rationality and $C_2$-cofiniteness of $C(L_{so(2n+1)}(k+2,0),L_{so(2n+1)}(k,0)\otimes L_{so(2n+1)}(2,0))$.

We first recall  some results in \cite{AP}.  Let $\mathbf A$ be a vector space spanned by $\Psi^{\pm}_{i}$,  $\Phi^{\pm}_{i}$, $\Psi_{2n+1}$, $\Phi_{2n+1}$, $i=1, \cdots, n$, and $(,)$ be the nondegenerate symmetric bilinear form on $\mathbf A$ such that
$$(\Psi^{\pm}_i, \Psi^{\mp}_j)=(\Phi^{\pm}_i, \Phi^{\mp}_j)=\delta_{i,j}, ~~~(\Psi^{\pm}_i, \Psi^{\pm}_j)=(\Phi^{\pm}_i, \Phi^{\pm}_j)=0,\ \ \ (\Psi^{\pm}_i, \Phi^{\pm}_j)=(\Psi^{\mp}_i, \Phi^{\pm}_j)=0,$$
$$(\Psi_{2n+1}, \Psi_{2n+1})=(\Phi_{2n+1}, \Phi_{2n+1})=1,\ \ \  (\Psi_{2n+1}, \Phi_{2n+1})=0,$$
$$(\Psi_{2n+1}, \Psi^{\pm}_i)=(\Phi_{2n+1}, \Psi^{\pm}_i)=0,\ \ \  (\Psi_{2n+1}, \Phi^{\pm}_i)=(\Phi_{2n+1}, \Phi^{\pm}_i)=0,$$
for $i,j=1, \cdots, n$. Then we have a fermion vertex operator superalgebra $M_{\hat{\mathbf A}}(1, 0)$ associated to $\hat{\mathbf A}$ (see Subsection \ref{fermion}). In the following, we will use $F_{2n+1}$ to denote the fermion vertex operator superalgebra $M_{\hat{\mathbf A}}(1, 0)$.

Let $F^{\Psi}_{n+\frac{1}{2}}$ (resp. $F^{\Phi}_{n+\frac{1}{2}}$) be the subalgebra of $F_{2n+1}$ generated by $\Psi_i^{\pm}(-\frac{1}{2})\1$, $\Psi_{2n+1}(-\frac{1}{2})\1$ (resp. $\Phi_i^{\pm}(-\frac{1}{2})\1$, $\Phi_{2n+1}(-\frac{1}{2})\1$), $i=1, \cdots, n$. Then the following results have been obtained in  \cite{FF} (see also Theorem 5 of \cite{AP}).
\begin{theorem}\label{FDb}
Let $n$ be a positive integer such that $n\geq 2$. Then\\
 (i) The vertex operator algebras $(F_{n+\frac{1}{2}}^{\Psi})^{even}$, $(F_{n+\frac{1}{2}}^{\Phi})^{even}$ are isomorphic to $L_{so(2n+1)}(1,0)$.
 (ii) $(F_{n+\frac{1}{2}}^{\Psi})^{odd}$ and $(F_{n+\frac{1}{2}}^{\Phi})^{odd}$ viewed as modules of $L_{so(2n+1)}(1,0)$ are isomorphic to $L_{so(2n+1)}(1,\Lambda_1)$.
\end{theorem}

Define the lattice
$$R_{2n+1}=\Z x_1+\cdots+\Z x_n+\Z y_1+\cdots+\Z y_n+\Z x,$$
$$(x_i, x_j)=(y_i, y_j)=\delta_{i,j},~~~(x_i, y_j)=0,$$
$$(x, x_j)=(x, y_j)=0,~~~(x, x)=1,$$
where $i, j=1, \cdots, n$. We set $z_{2k-1}=x_k, z_{2k}=y_k, z_{2n+1}=x$, where $k=1, \cdots, n$.

 Following \cite{AP},  we choose the 2-cocycle $\epsilon: R_{2n+1}\times R_{2n+1}\to \{\pm1\}$ determined by
 \begin{align}
 \epsilon(z_i, z_j)=\left\{
\begin{array}{ll}
1,& \text{~if~}i\leq j,\\
-1,& \text{~if~} i>j.
\end{array}
\right.
 \end{align}
 Then we have a lattice vertex operator superalgebra $V_{R_{2n+1}}$  associated to $R_{2n+1}$ (see Subsection \ref{lattice}). We shall need the following (non-standard) version of the boson-fermion correspondence \cite{Fe}, \cite{F}, \cite{K2}  (see also Theorem 6 of \cite{AP}).
 \begin{theorem}\label{bfb}
 Let $n$ be a positive integer. Then there exists a vertex superalgebra isomorphism $\phi_{2n+1}: F_{2n+1}\to V_{R_{2n+1}}$ such that
 \begin{align*}
 &\Psi^+_k(-\frac{1}{2})\1\mapsto \frac{1}{\sqrt{2}}(e^{x_k}+e^{-y_k}),\ \ \
 \Psi^-_k(-\frac{1}{2})\1\mapsto \frac{1}{\sqrt{2}}(e^{-x_k}+e^{y_k}),\\
 &\Phi^+_k(-\frac{1}{2})\1\mapsto \frac{i}{\sqrt{2}}(e^{x_k}-e^{-y_k}),\ \ \
 \Phi^-_k(-\frac{1}{2})\1\mapsto \frac{-i}{\sqrt{2}}(e^{-x_k}-e^{y_k}),\\
 &\Psi_{2n+1}(-\frac{1}{2})\1\mapsto \frac{1}{\sqrt{2}}(e^{x}+e^{-x}),\ \ \
  \Phi_{2n+1}(-\frac{1}{2})\1\mapsto \frac{i}{\sqrt{2}}(e^{x}-e^{-x}),
 \end{align*}
 for $k=1,\cdots, n$.
 \end{theorem}

 Based on Theorems \ref{FDb} and \ref{bfb}, the following results have been obtained in Theorem 7 of \cite{AP}.
 \begin{theorem}\label{subB}
 Let $n$ be a positive integer such that $n\geq 2$. Then\\
 (i) The subalgebra of $V_{R_{2n+1}}$ generated by
 \begin{align}\label{generatorb1}
 e^\gamma \text{~ and~} e^{-\gamma},
  \end{align}where $\gamma=x_1+\cdots+x_n+y_1+\cdots+y_n+x$, is isomorphic to the rank-one lattice vertex operator algebra $V_{\mathbf L}$ such that
   $$\mathbf L=\Z\gamma\text{~ and~ }(\gamma, \gamma)=2n+1.$$
 (ii)  The subalgebra of $V_{R_{2n+1}}$ generated by
 \begin{align}\label{generatorb2}
 &e^{x_k-y_l}+e^{x_l-y_k},~~~e^{x_k-x_l}+e^{y_l-y_k}, ~~~e^{-x_k+y_l}+e^{-x_l+y_k},\notag\\
 &e^{x-y_k}+e^{x_k-x},~~~e^{x-x_k}+e^{y_k-x},
 \end{align}
 for $k, l=1, \cdots, n$ and $k\neq l$, is isomorphic to $L_{so(2n+1)}(2, 0)$.\\
 (iii) The vertex operator algebra $V_{\mathbf L}\otimes L_{so(2n+1)}(2, 0)$ is isomorphic to the subalgebra of $V_{R_{2n+1}}$ generated by elements (\ref{generatorb1}) and (\ref{generatorb2}).
 \end{theorem}

 Furthermore, it was noticed in \cite{AP} that the root lattice of the Lie algebra $sl_{2n+1}$ can be realized as a sublattice of $R_{2n+1}$:
 \begin{align*}
 A_{2n}=&\Z(x_1-x_2)+\cdots+\Z(x_{n-1}-x_n)+\Z(x_n-x)+\Z(x-y_n)\\
 &+\Z(y_n-y_{n-1})+\cdots+\Z(y_2-y_1).
 \end{align*}
Let $\lambda_0=0$ and
\begin{align*}
&\lambda_i=x_1+\cdots+ x_i,~~~\text{~for~} 1\leq i\leq n,\\
&\lambda_{n+1}=x_1+\cdots+ x_n+x,\\
&\lambda_{n+i}=x_1+\cdots+ x_n+x+y_n+\cdots+ y_{n-i+2},~~~\text{~for~} 2\leq i\leq n.
\end{align*}
Then the following results have been obtained in \cite{AP}.
\begin{theorem}\label{decompb}
Let $n$ be a positive integer such that $n\geq 2$. Then\\
(i) $A_{2n}={\mathbf L}^{\perp}$.\\
(ii) $V_{R_{2n+1}}$ viewed as a module of $V_{\mathbf L}\otimes L_{sl_{2n+1}}(1, 0)$ has the following decomposition
\begin{align}
V_{R_{2n+1}}\cong\bigoplus_{i=0}^{2n}V_{\lambda_i+\mathbf{L}+A_{2n}}\cong \bigoplus_{i=0}^{2n} V_{\mathbf{L}+\frac{i}{2n+1}\gamma}\otimes L_{sl_{2n+1}}(1, \Lambda_i).
\end{align}
\end{theorem}

By Theorems \ref{subB}, \ref{decompb}, the affine vertex operator algebra $L_{so(2n+1)}(2, 0)$ is a subalgebra of $L_{sl_{2n+1}}(1, 0)$ (see Remark 4 of \cite{AP}). Moreover, the following result has been obtained in \cite{KW} (see also Example 7.3.4 of \cite{W}).
\begin{proposition}\label{decompwakib}
Let $n$ be a positive integer such that $n\geq 2$. Then $L_{sl_{2n+1}}(1, 0)$ viewed as a module of $L_{so(2n+1)}(2, 0)$ has the following decomposition
$$L_{sl_{2n+1}}(1, 0)\cong L_{so(2n+1)}(2, 0)\oplus L_{so(2n+1)}(2, 2\Lambda_1).$$
\end{proposition}

Combining Theorems \ref{FDb}, \ref{decompb}, we have
\begin{proposition}\label{evenb}
Let $n$ be a positive integer such that $n\geq 2$. Then\\
(i) The even part $F_{2n+1}^{even}$ of $F_{2n+1}$ contains a subalgebra  isomorphic to $L_{so(2n+1)}(1, 0)\otimes L_{so(2n+1)}(1, 0)$. Moreover, $F_{2n+1}^{even}$ viewed as a module of  $L_{so(2n+1)}(1, 0)\otimes L_{so(2n+1)}(1, 0)$ is isomorphic to $L_{so(2n+1)}(1, 0)\otimes L_{so(2n+1)}(1, 0)\oplus L_{so(2n+1)}(1, \Lambda_1)\otimes L_{so(2n+1)}(1, \Lambda_1)$.\\
(ii)  The even part $V_{R_{2n+1}}^{even}$ of $V_{R_{2n+1}}$ contains the subalgebra $V_{2\mathbf{L}}\otimes L_{sl_{2n+1}}(1, 0)$. Moreover, $V_{R_{2n+1}}^{even}$ viewed as a module of $V_{2\mathbf L}\otimes L_{sl_{2n+1}}(1, 0)$ has the following decomposition
\begin{align*}
V_{R_{2n+1}}^{even}&\cong \bigoplus_{\substack{0\leq i\leq n;\\i: \text{even}}}V_{2\mathbf{L}+\frac{i}{2n+1}\gamma}\otimes L_{sl_{2n+1}}(1, \Lambda_i)\bigoplus \bigoplus_{\substack{1\leq i\leq n;\\n+i: \text{even}}}V_{2\mathbf{L}+\frac{n+i}{2n+1}\gamma}\otimes L_{sl_{2n+1}}(1, \Lambda_{n+i})\\
& \bigoplus\bigoplus_{\substack{0\leq i\leq n;\\i: \text{odd}}}V_{2\mathbf{L}+\frac{2n+1+i}{2n+1}\gamma}\otimes L_{sl_{2n+1}}(1, \Lambda_i)\bigoplus \bigoplus_{\substack{1\leq i\leq n;\\n+i: \text{odd}}}V_{2\mathbf{L}+\frac{3n+1+i}{2n+1}\gamma}\otimes L_{sl_{2n+1}}(1, \Lambda_{n+i}).
\end{align*}
\end{proposition}

We are now ready to determine  the commutant of   $L_{so(2n+1)}(2,0)$ in  $L_{so(2n+1)}(1,0)\otimes L_{so(2n+1)}(1,0)\oplus L_{so(2n+1)}(1,\Lambda_1)\otimes L_{so(2n+1)}(1,\Lambda_1)$. 
\begin{theorem}
Let $n$ be a positive integer such that $n\geq 2$.
Then the commutant of   $L_{so(2n+1)}(2,0)$ in  $L_{so(2n+1)}(1,0)\otimes L_{so(2n+1)}(1,0)\oplus L_{so(2n+1)}(1,\Lambda_1)\otimes L_{so(2n+1)}(1,\Lambda_1)$ is isomorphic to the lattice vertex algebra $V_{2\mathbf{L}}$.
\end{theorem}
\pf By the similar argument as in Theorem \ref{cosetD}, the assertion follows from Theorem \ref{bfb} and Propositions \ref{decompwakib}, \ref{evenb}.
\qed
\subsection{Rationality and $C_2$-cofiniteness of diagonal coset vertex operator algebras $C(L_{so(2n+1)}(k+2,0),L_{so(2n+1)}(k,0)\otimes L_{so(2n+1)}(2,0))$}
In this subsection, we study rationality and $C_2$-cofiniteness of $C(L_{so(2n+1)}(k+2,0),L_{so(2n+1)}(k,0)\otimes L_{so(2n+1)}(2,0)).$
\begin{theorem}\label{rationalb}
 Let $n, k$ be positive integers such that $n\geq 2$. Suppose that $$C(L_{so(2n+1)}(k+1,0),L_{so(2n+1)}(k,0)\otimes L_{so(2n+1)}(1,0)),$$ $$C(L_{so(2n+1)}(k+2,0),L_{so(2n+1)}(k+1,0)\otimes L_{so(2n+1)}(1,0))$$ are rational and $C_2$-cofinite. Then  $C(L_{so(2n+1)}(k+2,0),L_{so(2n+1)}(k,0)\otimes L_{so(2n+1)}(2,0))$ is rational and $C_2$-cofinite.
\end{theorem}
\pf Consider the vertex operator algebra $L_{so(2n+1)}(k,0)\otimes F_{2n+1}^{even}$. By Proposition \ref{evenb}, $L_{so(2n+1)}(k,0)\otimes F_{2n+1}^{even}$ contains a subalgebra $$L_{so(2n+1)}(k,0)\otimes L_{so(2n+1)}(1,0)\otimes L_{so(2n+1)}(1,0).$$ Moreover, viewed as an $L_{so(2n+1)}(k,0)\otimes L_{so(2n+1)}(1,0)\otimes L_{so(2n+1)}(1,0)$-module, the vertex operator algebra $L_{so(2n+1)}(k,0)\otimes F_{2n+1}^{even}$ has the following decomposition
 \begin{align*}
 L_{so(2n+1)}(k,0)\otimes &L_{so(2n+1)}(1,0)\otimes L_{so(2n+1)}(1,0)\\
 &\oplus L_{so(2n+1)}(k,0)\otimes L_{so(2n+1)}(1,\Lambda_1)\otimes L_{so(2n+1)}(1,\Lambda_1).
 \end{align*}
Therefore, $L_{so(2n+1)}(k,0)\otimes F_{2n+1}^{even}$ is an extension of
\begin{align*}
L_{so(2n+1)}(k+2,0)&\otimes C(L_{so(2n+1)}(k+2,0), L_{so(2n+1)}(k+1,0)\otimes L_{so(2n+1)}(1,0))\\
&\otimes C(L_{so(2n+1)}(k+1,0), L_{so(2n+1)}(k,0)\otimes L_{so(2n+1)}(1,0)).
 \end{align*}
 As a result,  $C(L_{so(2n+1)}(k+2,0), L_{so(2n+1)}(k,0)\otimes F_{2n+1}^{even})$ is an extension of
\begin{align*}
C(L_{so(2n+1)}(k+2,0), &L_{so(2n+1)}(k+1,0)\otimes L_{so(2n+1)}(1,0))\\
&\otimes C(L_{so(2n+1)}(k+1,0), L_{so(2n+1)}(k,0)\otimes L_{so(2n+1)}(1,0)).
 \end{align*}
  By the assumption and Theorems  \ref{abd}, \ref{extensionr}, $C(L_{so(2n+1)}(k+2,0), L_{so(2n+1)}(k,0)\otimes F_{2n+1}^{even})$ is rational and $C_2$-cofinite. By Theorem \ref{bfb}, $F_{2n+1}^{even}$ is isomorphic to $V_{R_{2n+1}}^{even}$. Hence, $C(L_{so(2n+1)}(k+2,0), L_{so(2n+1)}(k,0)\otimes V_{R_{2n+1}}^{even})$ is rational and $C_2$-cofinite. Furthermore,  we could show that $C(L_{so(2n+1)}(k+2,0), L_{so(2n+1)}(k,0)\otimes V_{R_{2n+1}}^{even})$ is strongly regular by the similar argument as in Theorem \ref{rationalD}.

On the other hand, by Propositions \ref{decompwakib}, \ref{evenb},  $L_{so(2n+1)}(k,0)\otimes V_{R_{2n+1}}^{even}$ contains a subalgebra isomorphic to $L_{so(2n+1)}(k,0)\otimes L_{so(2n+1)}(2,0)\otimes V_{2\mathbf{L}}$. Hence, $L_{so(2n+1)}(k,0)\otimes V_{R_{2n+1}}^{even}$ is an extension of  $$L_{so(2n+1)}(k+2,0)\otimes C(L_{so(2n+1)}(k+2,0), L_{so(2n+1)}(k,0)\otimes L_{so(2n+1)}(2,0))\otimes V_{2\mathbf{L}}.$$ Since $C(L_{so(2n+1)}(k+2,0), L_{so(2n+1)}(k,0)\otimes V_{R_{2n+1}}^{even})$ is strongly regular, by Lemma \ref{cosetr}, we have
 \begin{align*}
 &C(V_{2\mathbf{L}}, C(L_{so(2n+1)}(k+2,0), L_{so(2n+1)}(k,0)\otimes V_{R_{2n+1}}^{even}))\otimes V_{2\mathbf{L}}\\
 &=C(L_{so(2n+1)}(k+2,0)\otimes V_{2\mathbf{L}}, L_{so(2n+1)}(k,0)\otimes V_{R_{2n+1}}^{even})\otimes V_{2\mathbf{L}}
  \end{align*}is rational and $C_2$-cofinite. Moreover, $$C(L_{so(2n+1)}(k+2,0)\otimes V_{2\mathbf{L}}, L_{so(2n+1)}(k,0)\otimes V_{R_{2n+1}}^{even})\otimes V_{2\mathbf{L}}$$  is of CFT type, it follows from Theorem 4.5 of \cite{ABD} that $$C(L_{so(2n+1)}(k+2,0)\otimes V_{2\mathbf{L}}, L_{so(2n+1)}(k,0)\otimes V_{R_{2n+1}}^{even})\otimes V_{2\mathbf{L}}$$  is regular. It is known \cite{DLM1} that $V_{2\mathbf{L}}$ is a regular vertex operator algebra. Hence, $C(L_{so(2n+1)}(k+2,0)\otimes V_{2\mathbf{L}}, L_{so(2n+1)}(k,0)\otimes V_{R_{2n+1}}^{even})$ is regular by Proposition \ref{regulart}.  Furthermore,  we could show that $C(L_{so(2n+1)}(k+2,0)\otimes V_{2\mathbf{L}}, L_{so(2n+1)}(k,0)\otimes V_{R_{2n+1}}^{even})$ is strongly regular by the similar argument as in Theorem \ref{rationalD}.

Note that $C(L_{so(2n+1)}(k+2,0)\otimes V_{2\mathbf{L}}, L_{so(2n+1)}(k,0)\otimes V_{R_{2n+1}}^{even})$ is an extension of $$C(L_{so(2n+1)}(k+2,0), L_{so(2n+1)}(k,0)\otimes L_{so(2n+1)}(2,0)).$$ We now determine the decomposition of  $$C(L_{so(2n+1)}(k+2,0)\otimes V_{2\mathbf{L}}, L_{so(2n+1)}(k,0)\otimes V_{R_{2n+1}}^{even})$$ as a $C(L_{so(2n+1)}(k+2,0), L_{so(2n+1)}(k,0)\otimes L_{so(2n+1)}(2,0))$-module. By Propositions \ref{decompwakib}, \ref{evenb}, $C(L_{so(2n+1)}(k+2,0)\otimes V_{2\mathbf{L}}, L_{so(2n+1)}(k,0)\otimes V_{R_{2n+1}}^{even})$ has the following decomposition
\begin{align}\label{decompcosetB}
&C(L_{so(2n+1)}(k+2,0)\otimes V_{2\mathbf{L}}, L_{so(2n+1)}(k,0)\otimes V_{R_{2n+1}}^{even})\notag\\
&=C(L_{so(2n+1)}(k+2,0), L_{so(2n+1)}(k,0)\otimes L_{so(2n+1)}(2,0))\notag\\
&\ \ \ \oplus C(L_{so(2n+1)}(k+2,0), L_{so(2n+1)}(k,0)\otimes L_{so(2n+1)}(2,2\Lambda_1)).
\end{align}

To show that $C(L_{so(2n+1)}(k+2,0), L_{so(2n+1)}(k,0)\otimes L_{so(2n+1)}(2,0))$ is rational and $C_2$-cofinite, we define a linear map $\tilde\sigma\in \End(L_{so(2n+1)}(k,0)\otimes L_{sl_{2n+1}}(1,0))$ such that
$$\tilde\sigma|_{L_{so(2n+1)}(k,0)\otimes L_{so(2n+1)}(2,0)}=\id, \ \ \  \tilde\sigma|_{L_{so(2n+1)}(k,0)\otimes L_{so(2n+1)}(2,2\Lambda_1)}=-\id.$$
Note that $L_{so(2n+1)}(2,2\Lambda_1)$ is a simple current $L_{so(2n+1)}(2,0)$-module (see Remark 5 of \cite{AP}). Then $\tilde\sigma$ is an automorphism of   $L_{so(2n+1)}(k,0)\otimes L_{sl_{2n+1}}(1,0)$. Denote  by $\sigma$ the restriction of $\tilde\sigma$ on \begin{align*}
&C(L_{so(2n+1)}(k+2,0), L_{so(2n+1)}(k,0)\otimes L_{sl_{2n+1}}(1,0))\\
&=C(L_{so(2n+1)}(k+2,0), L_{so(2n+1)}(k,0)\otimes (L_{so(2n+1)}(2,0)\oplus L_{so(2n+1)}(2,2\Lambda_1))).
 \end{align*}
 Then, by the formula (\ref{decompcosetB}), $\sigma$ is an automorphism of $$C(L_{so(2n+1)}(k+2,0)\otimes V_{2\mathbf{L}}, L_{so(2n+1)}(k,0)\otimes V_{R_{2n+1}}^{even})$$ such that
\begin{align*}
C(L_{so(2n+1)}(k+2,0)&\otimes V_{2\mathbf{L}}, L_{so(2n+1)}(k,0)\otimes V_{R_{2n+1}}^{even})^{\langle\sigma\rangle}\\
&=C(L_{so(2n+1)}(k+2,0), L_{so(2n+1)}(k,0)\otimes L_{so(2n+1)}(2,0)).
\end{align*}
We have proved that $C(L_{so(2n+1)}(k+2,0)\otimes V_{2\mathbf{L}}, L_{so(2n+1)}(k,0)\otimes V_{R_{2n+1}}^{even})$ is strongly regular. Hence, by Theorem \ref{orbifold}, $C(L_{so(2n+1)}(k+2,0), L_{so(2n+1)}(k,0)\otimes L_{so(2n+1)}(2,0))$ is rational and $C_2$-cofinite.
\qed

\begin{remark}
  For positive integers $k, l$,  it is conjectured in \cite{CL} that the vertex operator algebra  $C(L_{so(2n+1)}(k+l,0),L_{so(2n+1)}(k,0)\otimes L_{so(2n+1)}(l,0))$ is  rational and $C_2$-cofinite. By Theorem \ref{rationalb}, if the conjecture is true for $l=1$ then it is true for $l=2$.
\end{remark}

\section{Rationality and $C_2$-cofiniteness of  diagonal coset vertex operator algebras $C(L_{sl_2}(k+4,0),L_{sl_2}(k,0)\otimes L_{sl_2}(4,0))$}
\def\theequation{5.\arabic{equation}}
\setcounter{equation}{0}

In this section, we shall show that $C(L_{sl_2}(k+4,0),L_{sl_2}(k,0)\otimes L_{sl_2}(4,0))$ is rational and $C_2$-cofinite if $ k$ is an admissible number for $\hat{sl_2}$.
\subsection{The commutant of $L_{sl_2}(4,0)$ in $L_{sl_2}(2, 0)\otimes L_{sl_2}(2, 0)\oplus L_{sl_2}(2, 2\Lambda_1)\otimes L_{sl_2}(2, 2\Lambda_1)$}
In this subsection, we will show that the commutant of  $L_{sl_2}(4,0)$ in the vertex operator algebra  $L_{sl_2}(2, 0)\otimes L_{sl_2}(2, 0)\oplus L_{sl_2}(2, 2\Lambda_1)\otimes L_{sl_2}(2, 2\Lambda_1)$ is isomorphic to a rank-one lattice vertex operator algebra. This fact will play an important role in  our study of rationality and $C_2$-cofiniteness of $C(L_{sl_2}(k+4,0),L_{sl_2}(k,0)\otimes L_{sl_2}(4,0))$.

First, we recall from \cite{AP} some facts about the vertex operator superalgebras $F_{\frac{3}{2}}^{\Psi}$, $F_{\frac{3}{2}}^{\Phi}$ and $V_{R_3}$. We will need the following well-known results about  $F_{\frac{3}{2}}^{\Psi}$ and $F_{\frac{3}{2}}^{\Phi}$ (see Proposition 5 of \cite{AP}).
\begin{theorem}\label{FDa}
 (i) The vertex operator algebras $(F_{\frac{3}{2}}^{\Psi})^{even}$ and $(F_{\frac{3}{2}}^{\Phi})^{even}$ are isomorphic to $L_{sl_2}(2,0)$.\\
 (ii) $(F_{\frac{3}{2}}^{\Psi})^{odd}$ and $(F_{\frac{3}{2}}^{\Phi})^{odd}$ viewed as modules of $L_{sl_2}(2,0)$ are isomorphic to $L_{sl_2}(2,2\Lambda_1)$.
\end{theorem}

We next recall from \cite{AP} some facts about the lattice vertex operator superalgebra $V_{R_3}$. Recall that the lattice $R_3=\Z x_1+\Z y_1+\Z x$. Set $\mu=x_1+y_1+x$ and $$\mathbb{L}=\Z \mu.$$ Then $\mathbb{L}$ is a sublattice of $R_3$. Recall that the root lattice of $sl_3$ can be realized as a sublattice of $R_3$:
$$A_2=\Z(x_1-x)+\Z(x-y_1).$$
The following results have been obtained in \cite{AP}.
\begin{theorem}\label{decompa}
(i) $A_2=\mathbb{L}^{\perp}$.\\
(ii) The lattice vertex operator superalgebra $V_{R_3}$ contains the subalgebra $V_{\mathbb{L}}\otimes V_{A_2}$. Moreover, $V_{R_3}$ viewed as a $V_{\mathbb{L}}\otimes V_{A_2}$-module has the following decomposition
\begin{align*}
V_{R_3}&\cong V_{\mathbb{L}+A_2}\oplus V_{x_1+\mathbb{L}+A_2}\oplus V_{x_1+x+\mathbb{L}+A_2}\\
&\cong V_{\mathbb{L}}\otimes L_{sl_3}(1, 0)\oplus V_{\frac{\mu}{3}+\mathbb{L}}\otimes L_{sl_3}(1, \Lambda_1)\oplus V_{\frac{2\mu}{3}+\mathbb{L}}\otimes L_{sl_3}(1, \Lambda_2).
\end{align*}
\end{theorem}

By Theorems \ref{FDa} and \ref{decompa},  we have
\begin{proposition}\label{evena}
(i) The even part $F_{3}^{even}$ of $F_{3}$ contains a subalgebra  isomorphic to $L_{sl_2}(2, 0)\otimes L_{sl_2}(2, 0)$. Moreover, $F_{3}^{even}$ viewed as a module of  $L_{sl_2}(2, 0)\otimes L_{sl_2}(2, 0)$ is isomorphic to $L_{sl_2}(2, 0)\otimes L_{sl_2}(2, 0)\oplus L_{sl_2}(2, 2\Lambda_1)\otimes L_{sl_2}(2, 2\Lambda_1)$.\\
(ii)  The even part $V_{R_{3}}^{even}$ of $V_{R_{3}}$ contains the subalgebra $V_{2\mathbb{L}}\otimes L_{sl_{3}}(1, 0)$. Moreover, $V_{R_{3}}^{even}$ viewed as a module of $V_{2\mathbb{L}}\otimes L_{sl_{3}}(1, 0)$ has the following decomposition
\begin{align*}
V_{R_{3}}^{even}\cong V_{2\mathbb{L}}\otimes L_{sl_3}(1, 0)\oplus V_{\frac{4\mu}{3}+2\mathbb{L}}\otimes L_{sl_3}(1, \Lambda_1)\oplus V_{\frac{2\mu}{3}+2\mathbb{L}}\otimes L_{sl_3}(1, \Lambda_2).
\end{align*}
\end{proposition}

It is well-known that the vertex operator algebra $L_{sl_3}(1, 0)$ contains a subalgebra isomorphic to $L_{sl_2}(4, 0)$. Moreover, the decomposition of $L_{sl_3}(1, 0)$ as an $L_{sl_2}(4, 0)$-module is well-known (see Note 7.3.2 of \cite{W}).
\begin{proposition}\label{decompwakia}
The vertex operator algebra $L_{sl_3}(1, 0)$ contains a subalgebra isomorphic to $L_{sl_2}(4, 0)$. Moreover, $L_{sl_3}(1, 0)$ viewed as an $L_{sl_2}(4, 0)$-module has the following decomposition
$$L_{sl_3}(1, 0)=L_{sl_2}(4, 0)\oplus L_{sl_2}(4, 4\Lambda_1).$$
\end{proposition}

We are now ready to determine the commutant of $L_{sl_2}(4,0)$ in the vertex operator algebra  $L_{sl_2}(2, 0)\otimes L_{sl_2}(2, 0)\oplus L_{sl_2}(2, 2\Lambda_1)\otimes L_{sl_2}(2, 2\Lambda_1)$.
\begin{theorem}
The commutant of $L_{sl_2}(4,0)$ in the vertex operator algebra  $L_{sl_2}(2, 0)\otimes L_{sl_2}(2, 0)\oplus L_{sl_2}(2, 2\Lambda_1)\otimes L_{sl_2}(2, 2\Lambda_1)$ is isomorphic to the lattice vertex operator algebra $V_{2\mathbb{L}}$.
\end{theorem}
\pf By the similar argument as in Theorem \ref{cosetD}, the assertion follows from Theorem \ref{bfb} and Propositions \ref{evena}, \ref{decompwakia}.
\qed
\subsection{Rationality and $C_2$-cofiniteness of diagonal coset vertex operator algebras $C(L_{sl_2}(k+4,0),L_{sl_2}(k,0)\otimes L_{sl_2}(4,0))$} In this subsection, we will show that $C(L_{sl_2}(k+4,0),L_{sl_2}(k,0)\otimes L_{sl_2}(4,0))$ is rational and $C_2$-cofinite if $k$ is an admissible number for $\hat{sl_2}$. We will  need the following result, which has been proved in Theorem 6.3 of \cite{CFL}.
\begin{theorem}\label{rationala}
 Let $ k$ be an admissible number for $\hat{sl_2}$. Then the vertex operator algebra $C(L_{sl_2}(k+2,0),L_{sl_2}(k,0)\otimes L_{sl_2}(2,0))$ is strongly regular.
\end{theorem}

We are now ready to prove the main result in this section.
\begin{theorem}\label{rationalmain}
 Let $ k$ be an admissible number for $\hat{sl_2}$. Then the vertex operator algebra $C(L_{sl_2}(k+4,0),L_{sl_2}(k,0)\otimes L_{sl_2}(4,0))$ is rational and $C_2$-cofinite.
\end{theorem}
\pf Consider the vertex operator algebra $L_{sl_2}(k,0)\otimes F_{3}^{even}$. By Proposition \ref{evena}, $L_{sl_2}(k,0)\otimes F_{3}^{even}$ contains a subalgebra $L_{sl_2}(k,0)\otimes L_{sl_2}(2,0)\otimes L_{sl_2}(2,0).$ Moreover,  $L_{sl_2}(k,0)\otimes F_{3}^{even}$ viewed as a module of $L_{sl_2}(k,0)\otimes L_{sl_2}(2,0)\otimes L_{sl_2}(2,0)$ has the decomposition $ L_{sl_2}(k,0)\otimes L_{sl_2}(2,0)\otimes L_{sl_2}(2,0)\oplus L_{sl_2}(k,0)\otimes L_{sl_2}(2,2\Lambda_1)\otimes L_{sl_2}(2,2\Lambda_1).$
Therefore, $L_{sl_2}(k,0)\otimes F_{3}^{even}$ is an extension of
\begin{align*}
L_{sl_2}(k+4,0)&\otimes C(L_{sl_2}(k+4,0), L_{sl_2}(k+2,0)\otimes L_{sl_2}(2,0))\\
&\otimes C(L_{sl_2}(k+2,0), L_{sl_2}(k,0)\otimes L_{sl_2}(2,0)).
 \end{align*}
 As a result,  $C(L_{sl_2}(k+4,0), L_{sl_2}(k,0)\otimes F_{3}^{even})$ is an extension of
\begin{align*}
C(L_{sl_2}(k+4,0), L_{sl_2}(k+2,0)\otimes L_{sl_2}(2,0))\otimes C(L_{sl_2}(k+2,0), L_{sl_2}(k,0)\otimes L_{sl_2}(2,0)).
 \end{align*}
  By Theorems  \ref{rationala}, \ref{abd} and \ref{extensionr}, $C(L_{sl_2}(k+4,0), L_{sl_2}(k,0)\otimes F_{3}^{even})$ is rational and $C_2$-cofinite. By Theorem \ref{bfb}, $F_{3}^{even}$ is isomorphic to $V_{R_{3}}^{even}$. Hence, the vertex operator algebra $C(L_{sl_2}(k+4,0), L_{sl_2}(k,0)\otimes R_{3}^{even})$ is rational and $C_2$-cofinite. Furthermore,  we could show that $C(L_{sl_2}(k+4,0), L_{sl_2}(k,0)\otimes R_{3}^{even})$ is strongly regular by the similar argument as in Theorem \ref{rationalD}.

On the other hand, by Propositions \ref{decompwakia}, \ref{evena},  $L_{sl_2}(k,0)\otimes R_{3}^{even}$ contains a subalgebra isomorphic to $L_{sl_2}(k,0)\otimes L_{sl_2}(4,0)\otimes V_{2\mathbb{L}}$. Hence, $L_{sl_2}(k,0)\otimes R_{3}^{even}$ is an extension of  $$L_{sl_2}(k+4,0)\otimes C(L_{sl_2}(k+4,0), L_{sl_2}(k,0)\otimes L_{sl_2}(4,0))\otimes V_{2\mathbb{L}}.$$ Since $C(L_{sl_2}(k+4,0), L_{sl_2}(k,0)\otimes R_{3}^{even})$ is strongly regular, by Lemma \ref{cosetr}, we have
\begin{align*}
&C(V_{2\mathbb{L}}, C(L_{sl_2}(k+4,0), L_{sl_2}(k,0)\otimes R_{3}^{even}))\otimes V_{2\mathbb{L}}\\
&=C(L_{sl_2}(k+4,0)\otimes V_{2\mathbb{L}}, L_{sl_2}(k,0)\otimes R_{3}^{even})\otimes V_{2\mathbb{L}}
 \end{align*} is rational and $C_2$-cofinite. Moreover, $C(L_{sl_2}(k+4,0)\otimes V_{2\mathbb{L}}, L_{sl_2}(k,0)\otimes R_{3}^{even})\otimes V_{2\mathbb{L}}$ is of CFT type, it follows from Theorem 4.5 of \cite{ABD} that $$C(L_{sl_2}(k+4,0)\otimes V_{2\mathbb{L}}, L_{sl_2}(k,0)\otimes R_{3}^{even})\otimes V_{2\mathbb{L}}$$  is regular. It is known \cite{DLM1} that $V_{2\mathbb{L}}$ is a regular vertex operator algebra. Hence, $C(L_{sl_2}(k+4,0)\otimes V_{2\mathbb{L}}, L_{sl_2}(k,0)\otimes R_{3}^{even})$ is regular by Proposition \ref{regulart}.  Furthermore,  we could show that $C(L_{sl_2}(k+4,0)\otimes V_{2\mathbb{L}}, L_{sl_2}(k,0)\otimes R_{3}^{even})$ is strongly regular by the similar argument as in Theorem \ref{rationalD}.

Note that $C(L_{sl_2}(k+4,0)\otimes V_{2\mathbb{L}}, L_{sl_2}(k,0)\otimes R_{3}^{even})$ is an extension of $$C(L_{sl_2}(k+4,0), L_{sl_2}(k,0)\otimes L_{sl_2}(4,0)).$$ We now determine the decomposition of  $C(L_{sl_2}(k+4,0)\otimes V_{2\mathbb{L}}, L_{sl_2}(k,0)\otimes R_{3}^{even})$ as a $C(L_{sl_2}(k+4,0), L_{sl_2}(k,0)\otimes L_{sl_2}(4,0))$-module. By Propositions \ref{decompwakia}, \ref{evena}, the vertex operator algebra $C(L_{sl_2}(k+4,0)\otimes V_{2\mathbb{L}}, L_{sl_2}(k,0)\otimes R_{3}^{even})$ has the following decomposition
\begin{align}\label{decompcoseta}
&C(L_{sl_2}(k+4,0)\otimes V_{2\mathbb{L}}, L_{sl_2}(k,0)\otimes R_{3}^{even})\\
&=C(L_{sl_2}(k+4,0), L_{sl_2}(k,0)\otimes L_{sl_2}(4,0))\oplus C(L_{sl_2}(k+4,0), L_{sl_2}(k,0)\otimes L_{sl_2}(4,4\Lambda_1))\notag.
\end{align}

To show that $C(L_{sl_2}(k+4,0), L_{sl_2}(k,0)\otimes L_{sl_2}(4,0))$ is rational and $C_2$-cofinite, we define a linear map $\tilde\tau\in \End(L_{sl_2}(k,0)\otimes L_{sl_{3}}(1,0))$ such that
$$\tilde\tau|_{L_{sl_2}(k,0)\otimes L_{sl_2}(4,0)}=\id, \ \ \  \tilde\tau|_{L_{sl_2}(k,0)\otimes L_{sl_2}(4,4\Lambda_1)}=-\id.$$
Note that $L_{sl_2}(4,4\Lambda_1)$ is a simple current module of $L_{sl_2}(4,0)$ (see Proposition 2.20 of \cite{L3}). Then $\tilde\tau$ is an automorphism of the vertex operator algebra  $L_{sl_2}(k,0)\otimes L_{sl_{3}}(1,0)$. Denote by $\tau$ the restriction of $\tilde\tau$ on \begin{align*}
&C(L_{sl_2}(k+4,0), L_{sl_2}(k,0)\otimes L_{sl_{3}}(1,0))\\
&=C(L_{sl_2}(k+4,0), L_{sl_2}(k,0)\otimes L_{sl_2}(4,0))\oplus C(L_{sl_2}(k+4,0), L_{sl_2}(k,0)\otimes L_{sl_2}(4,4\Lambda_1)).
\end{align*} Then, by the formula (\ref{decompcoseta}), $\tau$ is an automorphism of $$C(L_{sl_2}(k+4,0)\otimes V_{2\mathbb{L}}, L_{sl_2}(k,0)\otimes R_{3}^{even})$$ such that
\begin{align*}
C(L_{sl_2}(k+4,0)\otimes V_{2\mathbb{L}}, L_{sl_2}(k,0)\otimes R_{3}^{even})^{\langle\tau\rangle}=C(L_{sl_2}(k+4,0), L_{sl_2}(k,0)\otimes L_{sl_2}(4,0)).
\end{align*}
We have proved that $C(L_{sl_2}(k+4,0)\otimes V_{2\mathbb{L}}, L_{sl_2}(k,0)\otimes R_{3}^{even})$ is strongly regular. Then it follows from Theorem \ref{orbifold} that $C(L_{sl_2}(k+4,0), L_{sl_2}(k,0)\otimes L_{sl_2}(4,0))$ is rational and $C_2$-cofinite.
\qed
\begin{remark}
In Conjecture 6.2 of \cite{CFL}, $C(L_{sl_2}(k+l,0), L_{sl_2}(k,0)\otimes L_{sl_2}(l,0))$ is conjectured to be rational and $C_2$-cofinite if $k$ is an admissible number for $\hat{sl_2}$ and $l$ is a positive integer. In case $l=1$, the conjecture has been confirmed in \cite{DMZ}, \cite{Wa}. In case $l=2$, the conjecture has been confirmed in  \cite{A}, \cite{CFL}.
\end{remark}

\section{Classification of irreducible modules of   diagonal coset vertex operator algebras $C(L_{sl_2}(k+4,0),L_{sl_2}(k,0)\otimes L_{sl_2}(4,0))$}
\def\theequation{6.\arabic{equation}}
\setcounter{equation}{0}
In this section, we shall classify irreducible modules of $C(L_{sl_2}(k+4,0),L_{sl_2}(k,0)\otimes L_{sl_2}(4,0))$  when $k$ is a positive odd integer.

\subsection{Modules of diagonal coset vertex operator algebras} Let $\g$ be a finite dimensional simple Lie algebra and  $k, l$ be positive integers.
 In this subsection, we recall from \cite{Lin} some facts about modules of the vertex operator algebra $C(L_{\g}(k+l,0),L_{\g}(k,0)\otimes L_{\g}(l,0))$. For $\dot\Lambda\in P_+^{k}, \ddot\Lambda\in P_+^{l}$, it follows from Theorem \ref{moduleaff} that $L_{\g}(k, \dot{\Lambda})$ and $L_{\g}(l, \ddot{\Lambda})$ are $L_{\g}(k, 0)$-module and $L_{\g}(l, 0)$-module, respectively. Then $L_{\g}(k, \dot{\Lambda})\otimes L_{\g}(l, \ddot{\Lambda})$ is an $L_{\g}(k, 0)\otimes L_{\g}(l, 0)$-module. As a consequence, $L_{\g}(k, \dot{\Lambda})\otimes L_{\g}(l, \ddot{\Lambda})$ may be viewed as an $L_{\g}(k+l, 0)$-module \cite{Lin1}. Since $L_{\g}(k+l, 0)$ is strongly regular, $L_{\g}(k, \dot{\Lambda})\otimes L_{\g}(l, \ddot{\Lambda})$ is completely reducible as an $L_{\g}(k+l, 0)$-module. For any $\Lambda\in P_+^{k+l}$, we define $$M_{\dot{\Lambda}, \ddot{\Lambda}}^{\Lambda}=\Hom_{L_{\g}(k+l, 0)} (L_{\g}(k+l, \Lambda), L_{\g}(k, \dot{\Lambda})\otimes L_{\g}(l, \ddot{\Lambda})).$$ Then it was proved in \cite{Lin} that $C(L_{\g}(k+l,0),L_{\g}(k,0)\otimes L_{\g}(l,0))=M_{0, 0}^0$ and $M_{\dot{\Lambda}, \ddot{\Lambda}}^{\Lambda}$ is a $C(L_{\g}(k+l,0),L_{\g}(k,0)\otimes L_{\g}(l,0))$-module. Moreover, we have
  the following results, which were essentially established in \cite{KW} (see also \cite{Lin}).
\begin{proposition}\label{decomp}
 Then for $\dot\Lambda\in P_+^{k}, \ddot\Lambda\in P_+^{l}, \Lambda\in P_+^{k+l}$, we have\\
(1) $M_{\dot{\Lambda}, \ddot{\Lambda}}^{\Lambda}\neq 0$ if and only if $\dot{\Lambda}+\ddot{\Lambda}-\Lambda\in Q$.\\
% see page 194 of \cite{KW}
(2) $L_{\g}(k, \dot{\Lambda})\otimes L_{\g}(l, \ddot{\Lambda})$ viewed as an $L_{\g}(k+l, 0)\otimes C(L_{\g}(k+l,0),L_{\g}(k,0)\otimes L_{\g}(l,0))$-module has the following decomposition
$$L_{\g}(k, \dot{\Lambda})\otimes L_{\g}(l, \ddot{\Lambda})=\oplus_{\Lambda \in P_+^{k+l}; \dot{\Lambda}+\ddot{\Lambda}-\Lambda\in Q}L_{\g}(k+l, \Lambda)\otimes M_{\dot{\Lambda}, \ddot{\Lambda}}^{\Lambda}.$$
\end{proposition}

We next show that there may be isomorphisms between $C(L_{\g}(k+l,0),L_{\g}(k,0)\otimes L_{\g}(l,0))$-modules $\{M_{\dot{\Lambda}, \ddot{\Lambda}}^{\Lambda}|\dot\Lambda\in P_+^{k}, \ddot\Lambda\in P_+^{l}, \Lambda\in P_+^{k+l}\}$.  For any $h\in \h$, set
$$\Delta(h, z)=z^{h(0)}\exp\left(\sum_{n=1}^{\infty}\frac{h(n)(-z)^{-n}}{-n}\right).$$
Then the following result has been established in Proposition 5.4 of \cite{L4}.
\begin{proposition}
Let $M$ be an irreducible $L_{\g}(k, 0)$-module and $h$ be an element of $\h$ such that $h(0)$ has only integral eigenvalues on $L_{\g}(k, 0)$. Set
$$(M^{(h)}, Y_{M^{(h)}}(\cdot, z))=(M, Y(\Delta(h, z)\cdot, z)).$$ Then $(M^{(h)}, Y_{M^{(h)}}(\cdot, z))$ is an irreducible  $L_{\g}(k, 0)$-module.
\end{proposition}
 Let $h^1, \cdots, h^n$ be elements in $ \h$ defined by $\alpha_i(h^j)=\delta_{i,j}$, $i, j=1, \cdots, n$. Set $$P^{\vee}=\Z h^1+\cdots+\Z h^n.$$ Then it is known \cite{L3} that $h(0)$ has only integral eigenvalues on $L_{\g}(k, 0)$ if and only if $h\in P^{\vee}$. In particular, for any irreducible $L_{\g}(k, 0)$-module $M$ and $h\in P^{\vee}$, $M^{(h)}$ is also an irreducible $L_{\g}(k, 0)$-module.  We now let $Q^{\vee}=\Z\alpha_1^{\vee}+\cdots+\Z\alpha_n^{\vee}$. Then it was proved in Proposition 2.25 of \cite{L3} that $M^{(h)}\cong M$ for any irreducible $L_{\g}(k, 0)$-module $M$ and $h\in Q^{\vee}$. Therefore, this induces  an action of $P^{\vee}/Q^{\vee}$ on the set $\{L_{\g}(k, \Lambda)|\Lambda \in P_+^k\}$ (see Proposition 2.24 of \cite{L3}). Moreover, the following result has been established in Theorem 2.26 of \cite{L3}.
 \begin{proposition}\label{iden1}
 Let $\theta=\sum_{i=1}^n a_i\alpha_i$, $a_i\in \Z_+$, be the highest root of $\g$ and set $J=\{i|a_i=1\}$. Then $P^{\vee}/Q^{\vee}=\{ 0+Q^{\vee}\}\cup\{h^i+Q^{\vee}|i\in J\}$.
 \end{proposition}

Thus, for any $h+Q^{\vee}\in P^{\vee}/Q^{\vee}$, there exists an element $\Lambda^{(h)}\in P_+^k$ such that  $L_{\g}(k, \Lambda)^{(h)}$ is isomorphic to $L_{\g}(k, \Lambda^{(h)})$.  As a result, we obtain the following isomorphisms between $C(L_{\g}(k+l,0),L_{\g}(k,0)\otimes L_{\g}(l,0))$-modules by using the operator $\Delta(\cdot, z)$ (see Corollary 4.4 of \cite{Lin}).
\begin{proposition}\label{iden}
For any  $h+Q^{\vee}\in P^{\vee}/Q^{\vee}$ and $\dot\Lambda\in P_+^{k}, \ddot\Lambda\in P_+^{l}, \Lambda\in P_+^{k+l}$, we have $M_{\dot{\Lambda}, \ddot{\Lambda}}^{\Lambda}\cong M_{\dot{\Lambda}^{(h)}, \ddot{\Lambda}^{(h)}}^{\Lambda^{(h)}}$ as $C(L_{\g}(k+l,0),L_{\g}(k,0)\otimes L_{\g}(l,0))$-modules.
\end{proposition}

We next consider the set $$\Omega=\{(\dot\Lambda, \ddot\Lambda, \Lambda)|\dot\Lambda\in P_+^{k}, \ddot\Lambda\in P_+^{l}, \Lambda\in P_+^{k+l}~\text{such that }~\dot{\Lambda}+\ddot{\Lambda}-\Lambda\in Q\}.$$
By Propositions \ref{decomp}, \ref{iden}, we have $\dot\Lambda^{(h)}+ \ddot\Lambda^{(h)}- \Lambda^{(h)}\in Q$ if  $\dot{\Lambda}+\ddot{\Lambda}-\Lambda\in Q$ and $h+Q^{\vee}\in P^{\vee}/Q^{\vee}$. Therefore, we have  an action of $P^{\vee}/Q^{\vee}$ on $\Omega$ defined as follows:
\begin{align*}
\pi: P^{\vee}/Q^{\vee}\times \Omega&\to \Omega\\
(h+Q^{\vee}, (\dot\Lambda, \ddot\Lambda, \Lambda))&\mapsto (\dot\Lambda^{(h)}, \ddot\Lambda^{(h)}, \Lambda^{(h)}).
\end{align*}
For any $(\dot\Lambda, \ddot\Lambda, \Lambda)\in \Omega$, we define the stabilizer $(P^{\vee}/Q^{\vee})_{(\dot\Lambda, \ddot\Lambda, \Lambda)}$ of $(\dot\Lambda, \ddot\Lambda, \Lambda)\in \Omega$ by $$(P^{\vee}/Q^{\vee})_{(\dot\Lambda, \ddot\Lambda, \Lambda)}=\{h+Q^{\vee}\in P^{\vee}/Q^{\vee}|(\dot\Lambda, \ddot\Lambda, \Lambda)=(\dot\Lambda^{(h)}, \ddot\Lambda^{(h)}, \Lambda^{(h)})\}.$$ For any  $(\dot\Lambda, \ddot\Lambda, \Lambda)\in \Omega$, we denote the orbit of $(\dot\Lambda, \ddot\Lambda, \Lambda)$ by $[\dot\Lambda, \ddot\Lambda, \Lambda]$. Then the following results have been obtained in Theorem 4.8 of \cite{Lin1}.

\begin{theorem}\label{classi}
Let $\g$ be a finite dimensional simple Lie algebra and  $k, l$ be positive integers. Suppose that the vertex operator algebra $C(L_{\g}(k+l,0),L_{\g}(k,0)\otimes L_{\g}(l,0))$ is rational, $C_2$-cofinite and that the stabilizer $(P^{\vee}/Q^{\vee})_{(\dot\Lambda, \ddot\Lambda, \Lambda)}$ of $(\dot\Lambda, \ddot\Lambda, \Lambda)\in \Omega$ is trivial for any  $(\dot\Lambda, \ddot\Lambda, \Lambda)\in \Omega$. Then \\
(1) For any $(\dot\Lambda, \ddot\Lambda, \Lambda)\in \Omega$, $M_{\dot{\Lambda}, \ddot{\Lambda}}^{\Lambda}$ is an irreducible $C(L_{\g}(k+l,0),L_{\g}(k,0)\otimes L_{\g}(l,0))$-module.\\
(2) $\{M_{\dot{\Lambda}, \ddot{\Lambda}}^{\Lambda}|[\dot\Lambda, \ddot\Lambda, \Lambda]\in \Omega/(P^{\vee}/Q^{\vee})\}$ is the complete list of inequivalent irreducible modules of  $C(L_{\g}(k+l,0),L_{\g}(k,0)\otimes L_{\g}(l,0))$.
\end{theorem}
\subsection{Classification of irreducible modules of  diagonal coset vertex operator algebras $C(L_{sl_2}(k+4,0),L_{sl_2}(k,0)\otimes L_{sl_2}(4,0))$} In this subsection, we will classify irreducible modules of   $C(L_{sl_2}(k+4,0),L_{sl_2}(k,0)\otimes L_{sl_2}(4,0))$ when $k$ is a positive odd integer. We use $\Delta=\{\pm \alpha\}$ to denote the root system of $sl_2$. Then we have $P^k_{+}=\{\frac{s\alpha}{2}|0\leq s\leq k\}$. By Proposition \ref{iden1}, we have $P^\vee/Q^\vee=\{Q^\vee, h^1+Q^\vee\}.$ Moreover, it has been proved in \cite{DR} that $L_{sl_2}(k, \frac{s\alpha}{2})^{(h^1)}\cong L_{sl_2}(k, \frac{(k-s)\alpha}{2})$. Therefore, we have $$\Omega=\{(i, j, k)|0\leq i\leq k, 0\leq j\leq 4, 0\leq s\leq k+4, i+j-s\in 2\Z\}.$$In the following, we use $M_{i, j}^{s}$ to denote the space $M_{\frac{i\alpha}{2}, \frac{j\alpha}{2}}^{\frac{s\alpha}{2}}$.

Since we have proved in Theorem \ref{rationalmain} that $C(L_{sl_2}(k+4,0),L_{sl_2}(k,0)\otimes L_{sl_2}(4,0))$ is rational and $C_2$-cofinite, we could apply Theorem \ref{classi} to classify irreducible modules of  $C(L_{sl_2}(k+4,0),L_{sl_2}(k,0)\otimes L_{sl_2}(4,0))$.
\begin{theorem}
 Let $k$ be a positive odd integer. Then for any $0\leq i\leq k$, $0\leq j\leq 4$ and $0\leq s\leq k+4$, we have \\
 (i) If the integer $i+j-s$ is even, then $M_{i, j}^{s}$ is an irreducible module of $C(L_{sl_2}(k+4,0),L_{sl_2}(k,0)\otimes L_{sl_2}(4,0))$.\\
 (ii) $M_{i, j}^{s}$ is isomorphic to $M_{k-i, 4-j}^{k+4-s}$.\\
 (iii) $\{M_{i, j}^{s}|[i,j,s]\in \Omega/(P^\vee/Q^\vee)\}$ is the complete list of inequivalent irreducible modules of  $C(L_{sl_2}(k+4,0),L_{sl_2}(k,0)\otimes L_{sl_2}(4,0))$.
\end{theorem}
\pf Assertion (i) follows from Proposition \ref{decomp} and Theorem \ref{classi}. Assertion (ii) follows from Proposition \ref{iden}. Assertion (iii) follows from Theorem \ref{classi}.
\qed


\begin{thebibliography}{ABCD}
\bibitem
%[ABD]
{ABD}
Abe, T.,  Buhl, G.,  Dong, C.: Rationality, regularity and $C_2$-cofiniteness. {\em Trans. Amer. Math. Soc.}  {\bf 356}, 3391-3402 (2004)

\bibitem{A}
 Adamovic, D.: Rationality of Neveu-Schwarz vertex operator superalgebras. {\em Internat. Math. Res. Notices}, 865-874 (1997)

\bibitem{AP}
Adamovic, D., Perse, O.: On coset vertex algebras with central charge 1. {\em Math. Commun.} {\bf 15}, 143-157  (2010)

\bibitem{ACKL}
Arakawa, T., Creutzig,  T., Kawasetsu,  K.,   Linshaw, A.: Orbifolds and cosets of minimal $W$-algebras. {\em Comm. Math. Phys.} {\bf 355},  339-372 (2017)

\bibitem{ACL}
Arakawa, T., Creutzig,  T.,    Linshaw, A.: W-algebras as coset vertex algebras. {\em Invent. Math.} {\bf 218}, 145-195 (2019)

\bibitem{ALY}
 Arakawa, T., Lam, C., Yamada, H.: Parafermion vertex operator algebras and W-algebras. {\em Trans. Amer. Math. Soc.} {\bf 371},  4277-4301  (2019)

\bibitem
%[B]
{B} Borcherds, R.: Vertex algebras, Kac-Moody algebras, and the Monster.
{\it Proc. Natl. Acad. Sci. USA} {\bf 83}, 3068-3071  (1986)

\bibitem{CM}
Carnahan, S., Miyamoto, M.: Regularity of fixed-point vertex operator subalgebras. {arXiv:1603.05645}.

\bibitem{CFL}
Creutzig, T., Feigin, B., Linshaw, A.: $N = 4$ superconformal algebras
and diagonal cosets. {\em Internat. Math. Res.
Notices} to appear. arXiv:1910.01228.
\bibitem{CL}
Creutzig, T.,  Linshaw, A.: Trialities of orthosymplectic $W$-algebras. arXiv:2102.10224.

\bibitem{DL}
Dong C.,  Lepowsky, J.: Generalized Vertex
Algebras and Relative Vertex Operators. Progress in Mathematics,
112. {\em Birkh\"{a}user, Boston, Inc., Boston, MA,} (1993)

\bibitem
%[DLM0]
{DLM0}
 Dong, C., Li, H., Mason, G.: Simple currents and extensions of vertex operator algebras. {\em Comm. Math. Phys.} {\bf 180}, 671-707  (1996)

\bibitem
%[DLM2]
{DLM1} Dong, C., Li, H., Mason, G.:
Regularity of rational vertex operator algebras. {\em  Adv. Math.} {\bf 132}, 148-166 (1997)

\bibitem
%[DLM3]
{DLM2} Dong, C., Li, H., Mason, G.: Modular invariance
of trace functions in orbifold theory and generalized moonshine. {\em
Comm. Math. Phys.} {\bf 214}, 1-56  (2000)




\bibitem
{DMZ}
Dong, C., Mason, G., Zhu, Y.: Discrete series of the Virasoro algebra and the moonshine module. {\em Algebraic groups and their generalizations: quantum and infinite-dimensional methods (University Park, PA, 1991)}, 295-316, Proc. Sympos. Pure Math., 56, Part 2, {\em Amer. Math. Soc., Providence, RI,} (1994)

\bibitem{DR}
 Dong, C.,  Ren, L.: Representations of the parafermion vertex operator algebras. {\em Adv. Math.} {\bf 315}, 88-101  (2017)

 \bibitem{FaZa}

 Zamolodchikov, A., Fateev, V.: Nonlocal (parafermion) currents in two-dimensional conformal quantum field theory and self-dual critical points in $Z_N$-symmetric statistical systems. {\em Soviet Phys. JETP} {\bf 62}, 215-225 (1985)
\bibitem{Fe}
 Feingold, A.: Constructions of vertex operator algebras. {\em Algebraic groups and their generalizations: quantum and infinite-dimensional methods (University Park, PA, 1991)}, 317-336, Proc. Sympos. Pure Math., 56, Part 2, {\em Amer. Math. Soc., Providence, RI,} (1994)

 \bibitem{FF}
  Feingold, A., Frenkel, I.: Classical affine algebras. {\em Adv. in Math.} {\bf 56}, 117-172  (1985)


\bibitem{F}
 Frenkel, I.: Two constructions of affine Lie algebra representations and boson-fermion correspondence in quantum field theory. {\em J. Functional Analysis} {\bf 44}, 259-327  (1981)


\bibitem
%[FHL]
{FHL} Frenkel, I., Huang, Y.,  Lepowsky, J.: On axiomatic
approaches to vertx operator algebras and modules. {\em Mem. Amer. Math. Soc.} {\bf  104} (1993)

\bibitem{FLM}
Frenkel, I., Lepowsky, J., Meurman, A.: Vertex operator algebras and the Monster. Pure and Applied Mathematics, 134. {\em Academic Press, Inc., Boston, MA}, (1988)
\bibitem
%[FZ]
{FZ} Frenkel, I.,  Zhu, Y.: Vertex operator algebras
associated to representations of affine and Virasoro algebra. {\em Duke.
Math. J.} {\bf 66}, 123-168 (1992)

\bibitem{GS}
Goddard, P., Schwimmer, A.: Unitary construction of extended conformal algebras. {\em Phys. Lett. B} {\bf 206},  62-70  (1988)

 \bibitem
%[HKL]
{HKL} Huang, Y.,  Kirillov, A. Jr.,  Lepowsky, J.: Braided tensor categories and extensions of vertex operator algebras. {\em Comm. Math. Phys.} {\bf 337},  1143-1159 (2015)


\bibitem{JL1}
Jiang, C., Lin, Z.: The commutant of $L_{\hat{sl}_2}(n, 0)$
in the vertex operator algebra $L_{\hat{sl}_2}(1, 0)^{\otimes n}$. {\em Adv. Math.} {\bf 301}, 227-257  (2016)

\bibitem{JL2}
Jiang, C., Lin, Z.: Tensor decomposition, parafermions, level-rank
duality, and reciprocity laws for vertex operator algebras. {\em Trans. Amer. Math. Soc.} to appear. arXiv:1406.4191.
\bibitem
%[K]
{K} Kac, V.: Infinite dimensional Lie algebras. Third edition. {\em Cambridge University Press, Cambridge} (1990)

\bibitem{K2}
Kac, V.: Vertex algebras for beginners. Second edition. University Lecture Series, 10. {\em American Mathematical Society, Providence, RI}, (1998)

\bibitem
{KW} Kac, V., Wakimoto, M.: Modular and conformal invariance constraints in representation theory of affine algebras. {\em Adv. Math.} {\bf 70}, 156-236  (1988)

 \bibitem{KW1}
 Kac, V., Wakimoto, M.: Classification of modular invariant representations of affine algebras. {\em Infinite-dimensional Lie algebras and groups (Luminy-Marseille, 1988)}, 138-177, Adv. Ser. Math. Phys., 7, {\em World Sci. Publ., Teaneck, NJ,} (1989)

\bibitem{KWang}
Kac, V., Wang, W.: Vertex operator superalgebras and their representations. {\em Mathematical aspects of conformal and topological field theories and quantum groups (South Hadley, MA, 1992)}, 161-191, Contemp. Math., 175, {\em Amer. Math. Soc., Providence, RI}, (1994)
\bibitem
%[LL]
{LL} Lepowsky, J.,  Li, H.: Introduction to vertex operator
algebras and their representations. Progress in Mathematics, 227. {\em Birkh$\ddot{a}$user Boston, Inc., Boston, MA,} (2004)

\bibitem
{Li}
Li, H.:  Some finiteness properties of regular vertex operator algebras. {\em J. Algebra} {\bf 212}, 495-514  (1999)

\bibitem
{L1}
Li, H.: Symmetric invariant bilinear forms on vertex operator algebras. {\em J. Pure Appl. Algebra} {\bf 96}, 279-297  (1994)

\bibitem{L2}
 Li, H.: Local systems of vertex operators, vertex superalgebras and modules. {\em J. Pure Appl. Algebra} {\bf 109}, 143-195  (1996)



\bibitem
%[L3]
{L3}
 Li, H.: Certain extensions of vertex operator algebras of affine type. {\em Comm. Math. Phys.} {\bf 217},  653-696  (2001)

 \bibitem
%[L1]
{L4}  Li, H.: Local systems of twisted vertex operators, vertex operator superalgebras and twisted modules, in: Moonshine, the Monster, and related topics. {\em Contemp. Math.} {\bf 193}, 203--236. American Mathematical Sociaty, Providence, RI,
   (1995)
 \bibitem{Lin}
 Lin, X.: Quantum dimensions and irreducible modules of some diagonal coset vertex operator algebras. {\em Lett. Math. Phys.} {\bf 110}, 1363-1380  (2020)
 \bibitem{Lin1}
 Lin, X.: Trace functions and fusion rules of diagonal coset vertex operator algebras. arXiv: 2104.06785.
 \bibitem{M}
Miyamoto, M.:
$C_2$-cofiniteness of cyclic-orbifold models. {\em Comm. Math. Phys.} {\bf 335}, 1279-1286  (2015)

 \bibitem{W}
  Wakimoto, M.: Lectures on infinite-dimensional Lie algebra. {\em World Scientific Publishing Co., Inc., River Edge, NJ}, (2001)

  \bibitem
%[W]
{Wa} Wang, W.: Rationality of Virasoro vertex operator algebras.
{\em Internat. Math. Res. Notices} {\bf  7}, 197-211  (1993)

\end{thebibliography}
\end{document}